\documentclass[smallextended,numbook,runningheads]{svjour3}     
\smartqed  
\usepackage{mathptmx}      
\usepackage{amsmath,amssymb,amsfonts,amsbsy,graphics} 
\usepackage{xcolor}

\usepackage{lipsum}
\usepackage{graphicx}
\usepackage{epstopdf}
\usepackage[ruled,vlined]{algorithm2e}
\usepackage{subfigure}

\newtheorem{prop}[theorem]{Proposition}

\newcounter{qcounter}

\begin{document}
	
	\title{WENO based adaptive image zooming algorithm
	}
	
	
	\author{ Bojan Crnkovi\'{c} \and Jerko \v{S}kifi\'{c} \and  Tina Bosner   
	}
	
	
	\institute{T. Bosner \at
		Department of Mathematics, Faculty of Science, University of Zagreb, Bijeni\v{c}ka cesta 30, 10000 Zagreb, Croatia \\
		\email{tinab@math.hr}           
		\and
		B. Crnkovi\'{c} \at
		Department of Mathematics, University of Rijeka, Radmile Matej\v{c}i\'{c} 2, 51000 Rijeka, Croatia \\
		\email{bojan.crnkovic@uniri.hr}
		\and
		J. \v{S}kifi\'{c} \at Department for Fluid Mechanics and Computational Engineering, Faculty of Engineering, University of Rijeka, Vukovarska 58, 51000 Rijeka, Croatia \\
		\email{jerko.skific@uniri.hr}
	}
	
	\date{Received: date / Accepted: date}

	\maketitle
	
	\begin{abstract}
Image zooming or upsampling is a widely used tool in image processing and an essential step in many algorithms. 
Upsampling increases the number of pixels and introduces new information into the image, which can lead to numerical effects such as ringing artifacts, aliasing effects and blurring of the image. 
In this paper, we propose an efficient polynomial interpolation algorithm based on the WENO algorithm for image upsampling that provides high accuracy in smooth regions, preserves edges, and reduces aliasing effects. Although, this is not the first application of WENO interpolation for image resampling, it is designed to have comparable complexity and memory load with better image quality than separable WENO algorithm.

We show that the algorithm performs equally well on smooth 2D functions, artificial pixel art and real digital images. 
Comparison with similar methods on test images shows good results on standard metrics and also provides visually satisfactory results. 
Moreover, the low complexity of the algorithm is ensured by a small local approximation stencil and the appropriate choice of smoothness indicators.
		
	\end{abstract}
	
	\section{Introduction}
Zooming a digital image increases the number of pixels of an original raster image. 
The pixel values of the new image must be determined based on the original image by an approximation algorithm. 
This process is very common and very important in many applications, from medical imaging to gaming or electronic publishing, denoising, antialiasing, satellite-image zooming and geometric transformation etc. 
Zooming algorithms are often executed in real time and therefore must be very efficient and not cause significant numerical artifacts.

\begin{figure}
		\centering
		\includegraphics[width=0.4\textwidth]{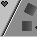}
		\includegraphics[width=0.4\textwidth]{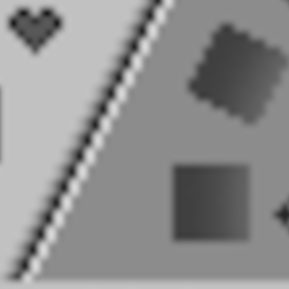}
		\caption{Artefacts and blurring of images}
		\label{figshadow}
	\end{figure}

We are interested in image upsampling efficient algorithm for image resolution doubling which preserves the edges and avoids ringing artifacts, aliasing effects, and excessive numerical diffusion. 
Ringing artifacts usually occur near sharp color transitions and look like ghost shadows of contours. 
This undesirable effect can be reduced by introducing artificial diffusion or by reducing oscillations in the approximation algorithm. 
Also, when lines are rendered in raster mode, a "staircase effect" can occur that can cause aliasing effects in the zoomed image as non-linear mixing effects create high-frequency components. 
Figure \ref{figshadow} shows a pixel graphics image and an 8x magnification with the standard bicubic algorithm, which exhibits all of the previously mentioned artifacts to some degree and introduces artificial diffusion

A common approach to image upsampling is to convert a discrete RGB image into a (often) continuous function with separate channels for each colour, and after reconstruction this function is transformed back into the new discrete image by resampling. So the main cause of ringing artefact is mostly oscillations of such approximations.	
There are a larger number of different approaches used for upsampling digital images \cite{Hou,Keys,Getreuer2011,Youngjoon2014,HakranKim2011} which can be divided into linear and nonlinear techniques.
The linear techniques, such as bilinear and bicubic spline interpolation \cite{Hou,Keys} have the advantage of simplicity and fast computation, but can lead to undesirable oscillations and/or reduced visual sharpness near sharp color transitions. On the other hand, non-linear methods can improve edges,
reduce artifacts and reconstruct pixel values with a high degree of accuracy, resulting in images with
higher subjective quality compared to linear methods. Although non-linear methods have obvious advantages over linear methods,
linear methods are used in most applications due to their low computational cost.

There are several good examples of non-linear methods based on the superposition of directed interpolation approximations.
Interpolation with geometric contour templates \cite{getreuergeomcontstenc} produces very sharp images with very good approximations of lines in multiple directions. This method reduces aliasing effects, but can produce strong ringing effects. Although this method produces visually sharp images, it is not suitable for all applications because it is very computationally intensive. The method in \cite{Lei2006} interpolates a missing sample in two orthogonal directions and then merges the results of the directional interpolation by the linear minimum mean square-error estimation.

Non-oscillatory methods developed primarily for solving PDEs have also been adapted for image denoising and resampling \cite{ImageZooming,Strong2000edge,Youngjoon2014,HakranKim2011}. These methods are computationally intensive but give very good results compared to the commonly used methods. Upsampling is treated as an inverse problem where the result is a plausible larger image that would look like the input image after downsampling. Good results are also obtained by \cite{Robidoux2008}, which leads to the histopolation problems.
The idea of histopolation is relatively old, and histopolation by polynomial splines was introduced by I. J. Schoenberg in \cite{schoenberghisto} in 1973. In \cite{bosnercrnkovicskifictension} a method based on reconstruction of surfaces using tension histosplines is presented. The basic idea of this method is to identify pixels with 2D numerical cells (instead of knots) and pixel values with cell averages. The histopolation approach to image reconstruction produces sharper images, but is more susceptible to numerical oscillations, so special treatment is required to reduce possible oscillations.
The tension splines we used in \cite{bosnercrnkovicskifictension} drastically reduce oscillations by applying tension to the splines, preserving the sharpness of the edges in the reconstructed surface that the image approximates.
The idea of splines with tension was also used in \cite{Yeon2015}, where the moving least squares method is controlled by a set of exponential polynomials with tension parameters so that they can be tuned to the characteristics of the given data. For a better fit to the local structures around the edges, the proposed algorithm also uses weights that take into account the edge orientation.

Some of 
the numerical methods developed for solving hyperbolic PDEs, such as  Weighted Essentially Non-Oscillatory (WENO), can be used for signal processing \cite{AMAT2021}  or resampling digital images \cite{ARANDIGA2012} where a non-separable two-dimensional weighted ENO interpolation was developed. The image can also be resampled using the tensor product of 1D WENO interpolation methods \cite{Mulet_2010}. These methods are highly accurate even in the presence of jump discontinuities and should not exhibit ringing effects. 

Many fast edge-enhancing methods (\cite{LocallyAdaptiveZooming,WenoZooming,Robidoux2009,Tian2012Edge}) for upsampling have two main stages: First, the main interpolation is used to create a double density version of the original image, and then this double density image is resampled using a simpler method to obtain an image with the desired resolution.

The main objective of this paper is to present a novel non-linear method optimized for image double density upsampling based on 1D WENO interpolation in multiple directions, resulting in a non-separable 2D approximation. The interpolation approximation in multiple directions are blended according to weights that depend on 1D smoothness indicators 
in the direction of interpolation. Although this is not the first use of WENO interpolation for digital image resampling \cite{ImageZooming,ARANDIGA2012,Mulet_2010}, this is the first one which is comparable with image quality and low computational complexity to methods of similar purpose available in standard libraries and open repositories.

The main advantage of this approach is to obtain a non-linear method with low complexity, able to adapt to local structures and to produce sharp images without oscillations and noticeable numerical artifacts. The complexity of the method is comparable to linear methods, but the resulting images have much better image quality. The method can be used for upsampling real digital images as well as for pixel art, and it can be used for interpolation of 2D functions evaluated on nested regular meshes. Although the algorithm works optimally for image doubling it can be used for arbitrary scale factors of magnification.

This paper is organized as follows. Section \ref{sec Interpolation} gives a brief introduction to fourth-order accurate 1D WENO interpolation, which serves as the basis and motivation for deriving 2D WENO as well for the tensor WENO interpolation. Section \ref{sed: 2D WENO} introduces a novel 2D WENO method with some theoretical results. Numerical experiments with artificial smooth 2D functions and the relative performance on real digital images is given in Section \ref{sec:numtests}, followed by a conclusion.

	\section{1D WENO interpolation 1D}
	\label{sec Interpolation}
	
	Let $a=x_{0}< \ldots <x_{N-1}=b$ be an equidistant partition of interval $[a,b]$, with step size $h$, and with given discrete set of function values $v_i=v(x_i),\: i=0, \ldots,N-1$  of some continuous function $v$. 
	WENO interpolation approximates $v$ by a rational function on $[x_i,x_{i+1}]$. 
	
	If a function $v$ is smooth enough, a polynomial interpolation can achieve high order of approximation on the interval $[x_i,x_{i+1}]$. 
	Let $p_i$ be an interpolating polynomial for nodes $\left\{x_{i-1},x_{i},x_{i+1}\right\}$ of degree  $2$:
	
	\begin{equation}\label{eq Tocnost}
	\begin{aligned}
		p_i(x)=v(x)+\mathcal{O}(h^{3}), \quad x \in [x_i,x_{i+1}],\\
		p_{i+1}(x)=v(x)+\mathcal{O}(h^{3}), \quad x \in [x_i,x_{i+1}].
	\end{aligned}
	\end{equation}
	 The following notation will be used for an interval: 
	\begin{equation}\label{eq Stencil}
		S_{i}=\left[x_{i-1},x_{i+1}\right]
	\end{equation}
	Both of intervals $S_i$ and $S_{i+1}$, contain $x_i$ and $x_{i+1}$, as can be seen in (\ref{eq Tocnost}), the approximations are of order $3$ on the interval $[x_i,x_{i+1}]$.
	The union of intervals $\mathcal{S}=S_i\cup S_{i+1}$, each contain $4$ nodes.
	 We can form an interpolating polynomial $q_i$ of degree $3$ on interval $\mathcal{S}_i$. If $v$ is smooth on interval $\mathcal{S}_i$,
	then $q_i$ is $\mathcal{O}(h^4)$ approximation of $v$. One can find polynomials  $C_{i,s},\ s=0,1$, so it holds:
	\begin{equation}\label{eq PolinomZaVelikiStencil}
		q(x)=
		C_{i,0}(x) p_i(x)+C_{i,1}(x) p_{i+1}(x), \quad x \in [x_i,x_{i+1}].
	\end{equation}
	$C_{i,s}(x)$ are first degree polynomials  
	which must satisfy consistency conditions: 
	$$C_{i,s}(x) \geq 0\textrm{ and } C_{i,0}(x)+C_{i,1}(x)=1, \quad  x \in [x_{i},x_{i+1}].$$
	From (\ref{eq PolinomZaVelikiStencil}) it follows that 
	we can approximate $v$ on $[x_i,x_{i+1}]$ with $4$\textsuperscript{th} order approximation as a convex combination of $3$\textsuperscript{rd} order accurate approximations. If $v$ is not smooth, loss of accuracy and oscillations in the approximation can appear.
	
	WENO interpolation can reduce oscillations of the approximation by measuring oscillations of polynomials $p_i$.
	Instead of  $q_i$ we will use a weighted combination of polynomials 
	
	\begin{equation}\label{eq TezinskiPolinom1d}
		u_i(x)=\omega_{i,0}(x) p_{i}(x)+ \omega_{i,1}(x) p_{i+1}(x),
	\end{equation}
	where $\omega_{i,s}$ are weighting functions which satisfy the same consistency conditions as ideal weights \eqref{eq PolinomZaVelikiStencil}. Furthermore, if the data is non-oscillatory on $S_i$, then $\omega_{i,s}$ must be very close to ideal weights  $C_{i,s}$:
	\begin{equation}\label{eq UvjetOmega-C_r-s}
		\omega_{i,s}(x)=C_{i,s}(x)+\mathcal{O}(h^{2}), \quad
		s=0,1.
	\end{equation}
	According to  \cite{shu_2020}, if the condition (\ref{eq UvjetOmega-C_r-s}) is satisfied, then
	$u_i$ is an approximation of order $4$. If the polynomial $p_i$ oscillates on 
	$[x_i,x_{i+1}]$ , then $\omega_{i,s}(x)$ should be very small to reduce oscillations in the final approximation. 
	
	Weights which satisfy given condition can be obtained in form of:
	\begin{equation}\label{eq IzrazOmega}
		\omega_{i,s}(x)=\frac{\alpha_{i,s}(x)}{\alpha_{i,0}(x)+\alpha_{i,1}(x)},\textrm{ and }
		\alpha_{i,s}(x)=\frac{C_{i,s}(x)}{\left(\varepsilon_h+SI_{i,s}\right)^{\beta}},\ s=0,1,
	\end{equation}
%
	where $\varepsilon_h=K_\varepsilon h^2$ prevents division by zero and is often very small, usually $K_\varepsilon=10^{-8}$, and $\beta>0$ amplifies how strongly interpolation reduces oscillations. 
	
	Smoothness indicator $SI_{i,s}$ which apers in (\ref{eq
		IzrazOmega}) measures oscillations of  $p_{i+s}$ on $[x_i,x_{i+1}]$.
	We use 
	the ones which satisfy (\ref{eq UvjetOmega-C_r-s}) :
	\begin{equation*}
		SI_{i,s}=\sum_{l=1}^{2} \int_{x_i}^{x_{i+1}} h^{2l-1}
		\left( \frac{\textrm{d}^{l}p_{i+s}(x)}{\textrm{d}x^{l}} \right)^{2} \,\textrm{d}x.
	\end{equation*}
	From the \cite{Liu2009} we have the following properties of 1D WENO smoothness indicators:
	\begin{equation}\label{eq: SIporperties}
		SI_{i,s}= 
		\begin{cases}
			\mathcal{O}\left(h^2\right)& \text{if } v \text{ is smooth locally,}\\
			\mathcal{O}\left(1\right),              & \text{otherwise.}
		\end{cases}
	\end{equation}
	More details about smoothness indicators can be found in \cite{shu_2020,Liu2009}.
	
	In special case of uniform mashes, when the approximation is evaluated in the middle of the interval $x=(x_i+x_{i+1})/2$, the ideal weight are equal $C_{i,0}=C_{i,1}=0.5$ which simplifies \eqref{eq IzrazOmega}.
	
	\section{{Weighted Direction WENO 2D}  interpolation }\label{sed: 2D WENO}
\begin{figure}
		\begin{center}
			\subfigure[Input data]{\includegraphics[width=0.43\textwidth]{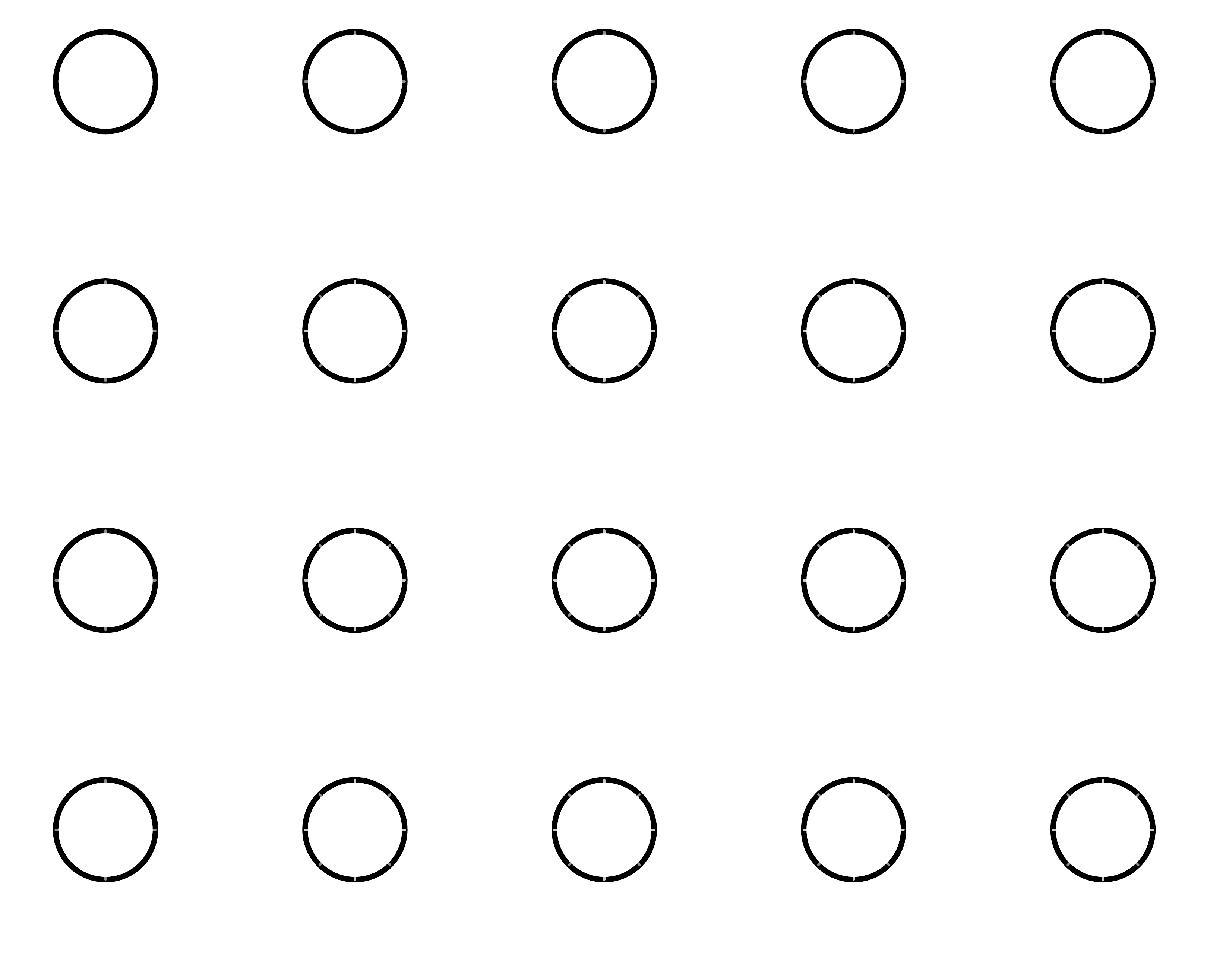}}
			\subfigure[Interpolated data]{\includegraphics[width=0.52\textwidth]{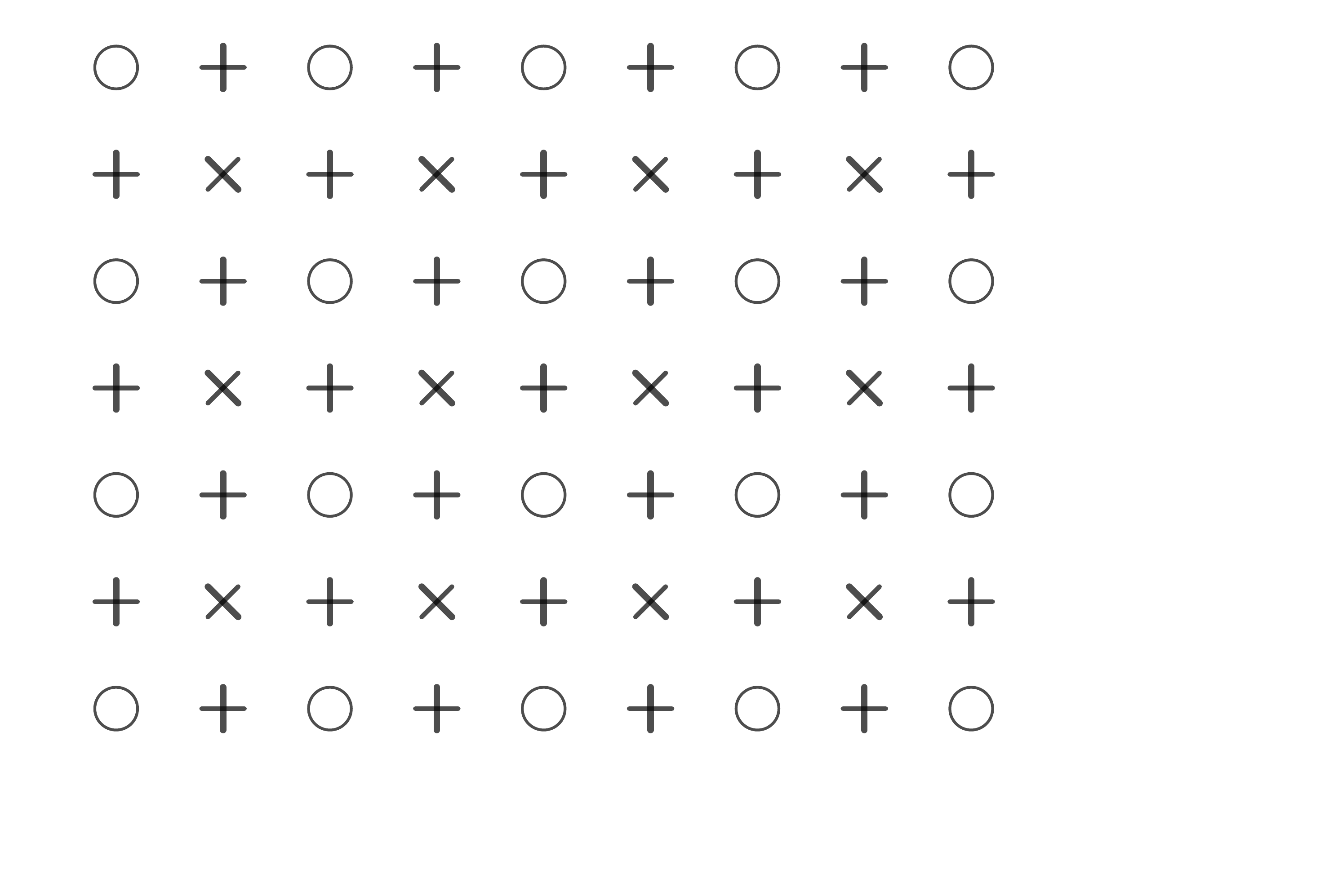}}
			\caption{Approximation of an $(2n-1)\times(2m-1)$ array from input $n\times m$ array }
			\label{fig:inputdata}
		\end{center}
	\end{figure}

	Using true 2D WENO interpolation can be very computationally demanding even on a uniform rectilinear meshes, so we will focus on a special case where the $n\times m$ rectilinear mesh, with step sizes in both dimensions equal to $h$, is interpolated to obtain a values on a finer mesh with $(2n-1)\times (2m-1)$ grid points. This will enable us to use fast 1D WENO interpolation from the previous section. Figure \ref{fig:inputdata} shows given points on a grid marked with circles and points which have to be interpolated marked with "$\times$" or "+". If we index given points with even coordinates $(2i,2j),\ 0\leq i<n,\ 0\leq j<m$, with the respect to the finer grid, then all "$\times$" marked points have odd coordinates $(2i+1,2j+1),\ 0\leq i<n-1,\ 0\leq j<m-1$ and all "+" marked points will have one even and the other odd coordinate.

	The interpolation is carried out in two phases. 
	The first phase carries out the interpolation in slanted directions (Figure \ref{fig:x}) and the second phase in horizontal and vertical directions (Figure \ref{fig:+}). Similarly to 1D case, if the function which we interpolate is smooth, the interpolation should be close to weighted combination of interpolating polynomials
	in directions 
	\begin{equation} \label{defgammak}
		\gamma_k=\gamma_0+k\pi/2, \quad k=0,1,2,3,
	\end{equation}
	which we denote
	\begin{equation}\label{eq: Ideal2DPolynomial}
		q(x_{i,j})=\delta\sum_{k=0}^1 C_{i,j}^{\gamma_{2k}} p_{i,j}^{\gamma_{2k}}+(1-\delta)\sum_{k=0}^1 C_{i,j}^{\gamma_{2k+1}} p_{i,j}^{\gamma_{2k+1}},
	\end{equation}
	where $x_{i,j}:=(x_i,y_j):=\left(x_0+i\frac{h}{2},y_0+j\frac{h}{2}\right)$, the weights are $C^{\gamma_{k}}_{i,j}=0.5$, $p_{i,j}^{\gamma_{k}}$ are quadratic interpolating polynomials in directions ${\gamma_{k}}$ evaluated at $x_{i,j}$, and $\delta$ is defined in (\ref{defdelta}). We should point out that the weighted combination of quadratic polynomials in parallel directions reconstruct a value of a cubic interpolating polynomial.
	
	\subsection{Phase 1}
	In the first phase, points with odd coordinates designated with "$\times$" are interpolated by a convex combination of values of second degree interpolating 1D polynomials in directions $\gamma=\frac{\pi}{4}, \frac{3\pi}{4}, \frac{5\pi}{4}, \frac{7\pi}{4}$, which correspond to (\ref{defgammak}) with $\gamma_0=\frac{\pi}{4}$. 
	The interpolating polynomials use given values marked with ("$\circ$") with even coordinates.

\begin{figure}
	\begin{center}
		\subfigure[Smoothness indicators associated with one direction]{\includegraphics[width=0.45\textwidth]{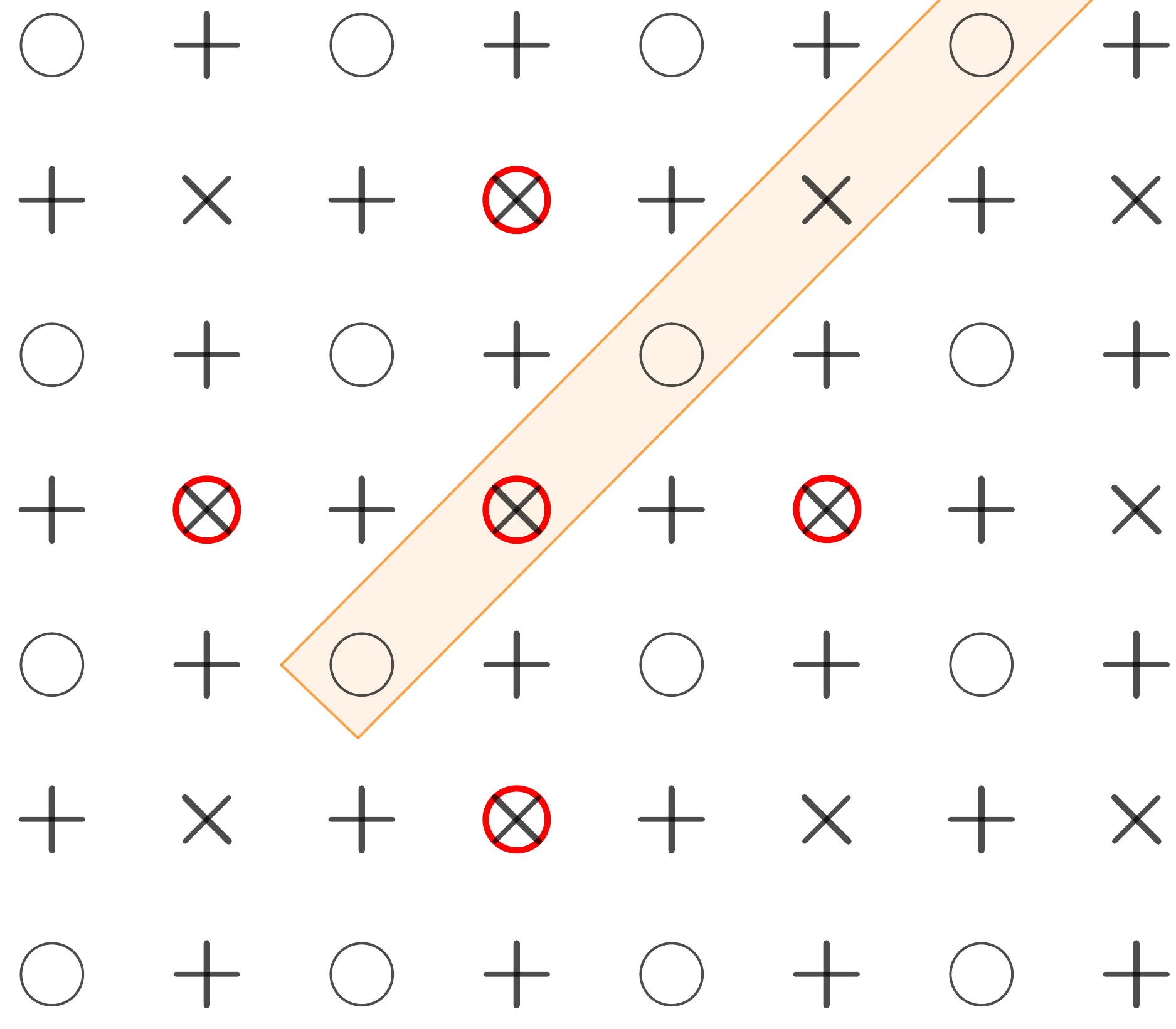}}
		\subfigure[Interpolation in slanted directions]{\includegraphics[width=0.45\textwidth]{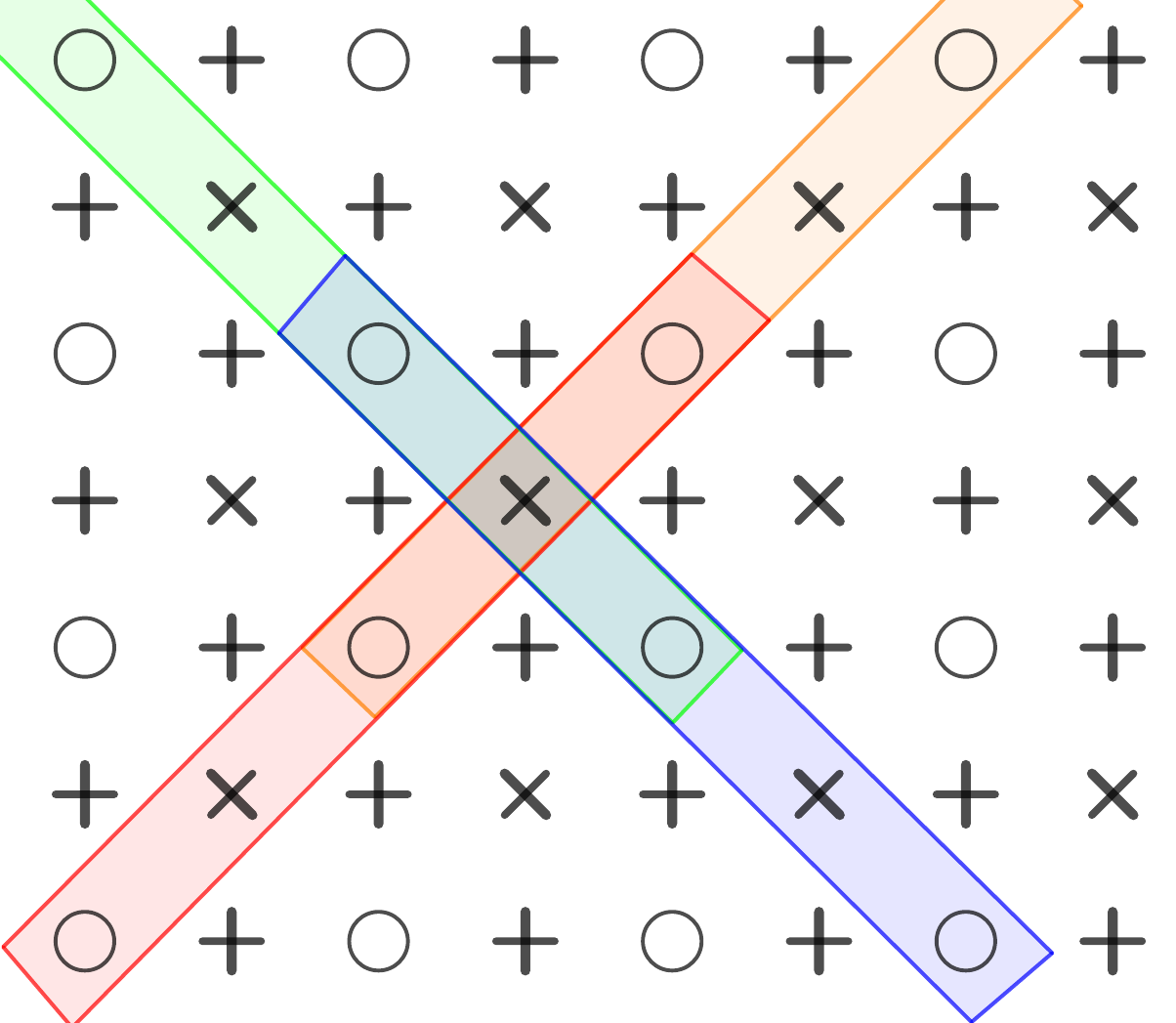}} \quad 
		\caption{Approximation of grid data marked with "$\times$" based on input data marked with "$\circ$".}
		\label{fig:x}
	\end{center}
\end{figure}	
	
	Figure \ref{fig:x} shows stencils for this 4 polynomials and we can clearly see that they overlap on the middle point marked with "$\times$".
	The idea behind this construction is that this interpolation should inherit good properties of WENO interpolation from 1D and if the data is not smooth it would pick the appropriate weight for the direction of interpolation in 2D. {This approach will be called Weighted Direction (WD) WENO interpolation further on.} 
	The key to this would be to pick correct weights for the linear combination of 1D interpolation values.
	
	Similar to the 1D idea, the values interpolated by four 1D polynomials are used to obtain a weighted combination for the point where all four polynomials intersect. We will introduce a vector which points in the direction of interpolation:
	\begin{equation}\label{eq:Vektor}
		d^{\gamma}=\frac{\sqrt{2}}{2}h\frac{1}{\cos^2(\gamma)+\sin^2(\gamma)}(\cos(\gamma),\sin(\gamma)).
	\end{equation}
	The interpolation stencil in 2D depends on direction of interpolation $\gamma$ and the central interpolation point
	\begin{equation}
		S_{i,j}^{\gamma}=\{x_{i,j}-d^{\gamma},x_{i,j}+d^{\gamma},x_{i,j}+3d^{\gamma}\}.
	\end{equation}
	We need to evaluate interpolation polynomials which depend on the stencil in the given direction
	\begin{equation}
		p_{i,j}^\gamma:=p_{S_{i,j}^{\gamma}}(x_{i,j}).
		\label{eq:p}
	\end{equation}
	We will then associate a 1D smoothness indicator to each polynomial which depends on the direction of interpolation: 
	\begin{equation}\label{eq 2Dsmoothness}
		SI_{i,j}^\gamma=\sum_{l=1}^{2} \int_{x_{i,j}-d^{\gamma}}^{x_{i,j}+d^{\gamma}} \Big(\sqrt{2}h\Big)^{2l-1}
		\left( \frac{\textrm{d}^{l}p_{S_{i,j}^{\gamma}}(x)}{\textrm{d}x^{l}} \right)^{2} \,\textrm{d}x{.}
	\end{equation}
	In 2D we will use a linear combination of nearby 1D smoothness indicators  \eqref{eq 2Dsmoothness} where 1D indicator  $SI_{i,j}^\gamma$ must remain if $h \to 0$:
	\begin{equation}\label{eq 2DsmoothnessFinal}
		D_{i,j}^\gamma= SI_{i,j}^\gamma+ \frac{h^2}{4}(SI_{i,j+2}^\gamma+SI_{i+2,j}^\gamma+SI_{i,j-2}^\gamma+SI_{i-2,j}^\gamma){.} 
	\end{equation}
	The equation \eqref{eq 2DsmoothnessFinal} shows that the smoothness indicator in {WD} WENO algorithm depends on smoothness of nearest polynomials calculated in this phase. Figure \ref{fig:+}a) shows in red the local dependence stencil of one of the 1D polynomials with smoothness indicators. This choice is not arbitrary because we need to make sure the {WD} WENO algorithm must have low computational load and give high order accurate approximations.
	  
	The ideal weights for this interpolation would be $0.5$ which would result in a 1D cubic interpolation polynomial. 
	However, we'd like to assign larger weight values to polynomials which have smaller values of smoothness indicator  \eqref{eq 2DsmoothnessFinal}:
	\begin{equation}
		\alpha_{i,j}^\gamma=\frac{0.5}{(\varepsilon_h+D_{i,j}^\gamma)^{\beta}}.
		\label{eq:alpha}
	\end{equation}
	
	Final non-linear weights used in the approximation must be scaled to sum up to one:
	\begin{equation}\label{eq IzrazOmega2d}
		\omega_{i,j}^{\gamma_k}=\frac{\alpha_{i,j}^{\gamma_k}}{\sum_{l=0}^3\alpha_{i,j}^{\gamma_l}}.
	\end{equation}
	
	Finally, interpolated value is obtained by a weighted combination of interpolation polynomials:
	 \begin{equation}\label{eq 2dKombinacija}
	 	u_{i,j}=\sum_{k=0}^3 \omega_{i,j}^{(2k+1)\pi/4} p_{i,j}^{(2k+1)\pi/4},
	 \end{equation}
	 where \begin{itemize}
	 	\item $\omega_{i,j}^\gamma$ are non-linear weights which depend on 1D smoothness indicators of nearby polynomials in the same direction, and
	 	\item $p_{i,j}^\gamma$ are associated values of the interpolation polynomials on stencil $S_{i,j}^{\gamma}$ at the intersecting point $x_{i,j}$.
	 \end{itemize}

\begin{algorithm}
	\SetAlgoLined
	\footnotesize
	\For {$k=0, \, k< 3,\, k=k+1 $}{
		$\gamma_k=\pi/4+k\pi/2$\\
		\For {$i=1, \,  i<2n-1, \,  i=i+2$}{
			\For {$j=1, \,  j<2m-1, \,  j=j+2$}{
				Calculate $p_{i,j}^{\gamma_k}$ Equation \eqref{eq:p}\\
				Calculate $SI_{i,j}^{\gamma_k}$ Equation \eqref{eq 2Dsmoothness}
			}
		}
		\For {$i=1, \,  i<2n-1, \,  i=i+2$}{
			\For {$j=1, \,  j<2m-1, \,  j=j+2$}{
				Calculate $\alpha_{i,j}^{\gamma_k}$ Equation (\ref{eq:alpha})\\
				$u_{i,j}=u_{i,j}+\alpha_{i,j}^{\gamma_k} p_{i,j}^{\gamma_k}$
			}
		}
	}
	\For {$i=1, \,  i<2n-1, \,  i=i+2$}{
		\For {$j=1, \,  j<2m-1, \,  j=j+2$}{
			$u_{i,j}=\frac{u_{i,j}}{\sum_{k=0}^3 \alpha_{i,j}^{\gamma_k}} $
		}
	}
	\caption{ {WD} WENO phase 1, $\times$ directions }
	\label{alg1}
\end{algorithm}

	\subsection{Phase 2}

The algorithm \ref{alg1} is carried out in the first phase where the "$\times$" marked points must be interpolated because they are needed for the second phase.

\begin{figure}
		\begin{center}
			{
				\subfigure[Smoothness indicator]{\includegraphics[width=0.42\textwidth]{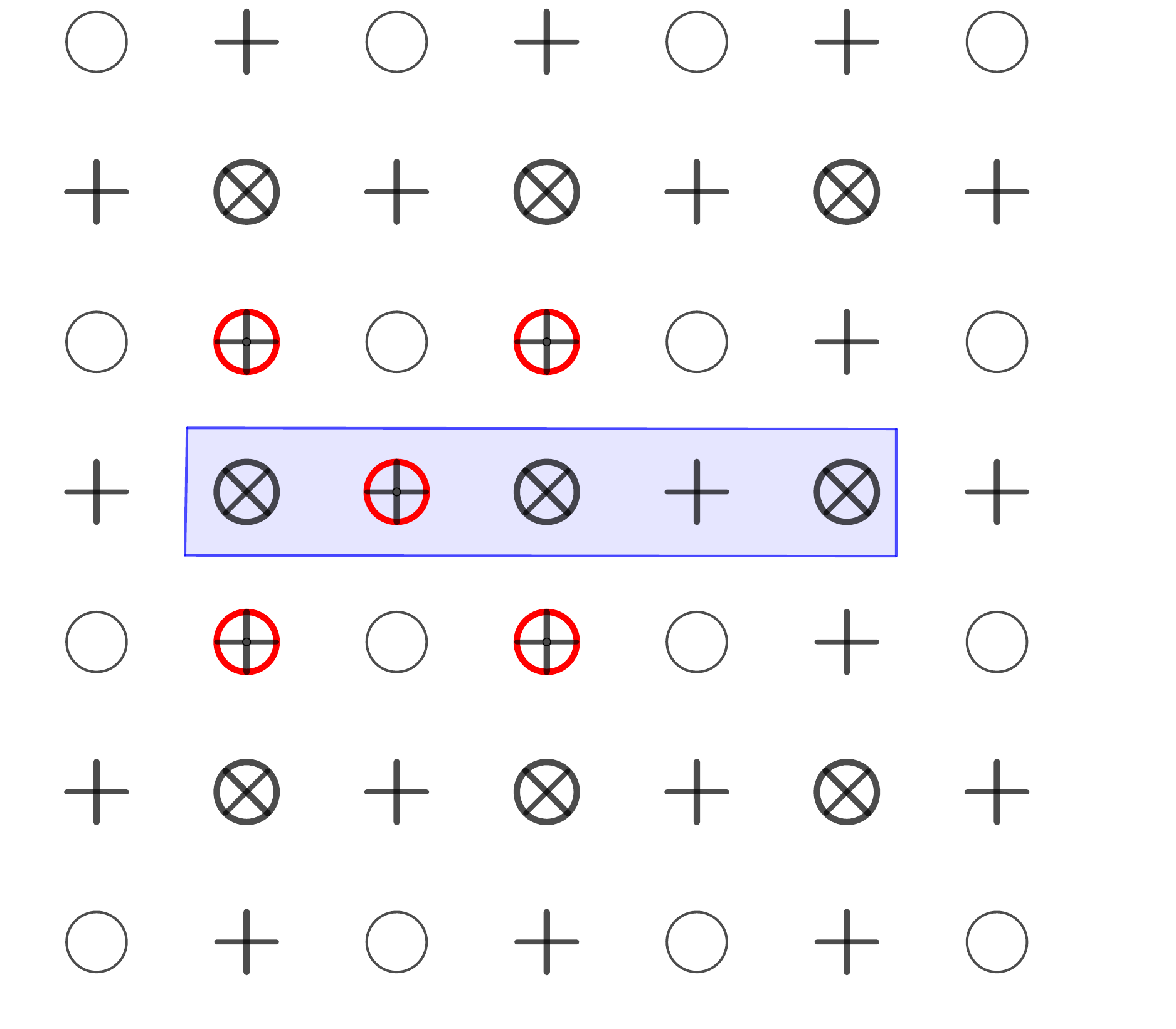}}
				\subfigure[Interpolation in horizontal and vertical directions]{\includegraphics[width=0.48\textwidth]{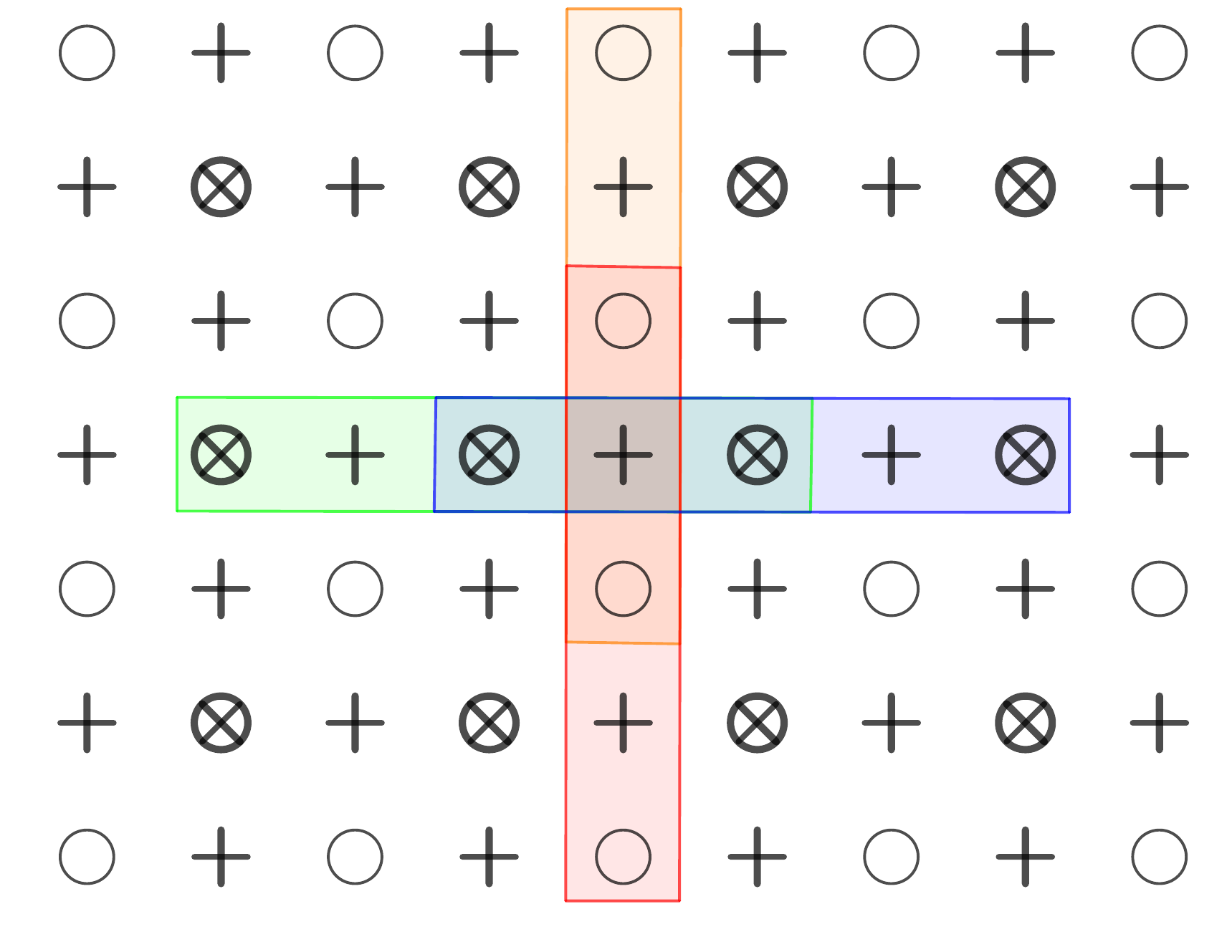}} }			
			\caption{Interpolation of grid data marked with "+".}
			\label{fig:+}
		\end{center}
\end{figure}

The second phase must interpolate "+" points and is carried out in similar to first phase. The only difference between these phases is the direction of interpolation, now $\gamma_0=0$ in (\ref{defgammak}). Figure \ref{fig:+} shows stencils for this 4 polynomials and we can clearly se that they overlap on the middle point marked with "+". Interpolation uses only the original local data on even coordinates and the data on points marked with "$\times$" on odd coordinates which were obtained in the previous step.
	
	Finally, interpolated value for coordinates with one odd and one even is obtained by a weighted combination of interpolation polynomials in standard coordinate directions where now we can use values obtained in previous stage:
	
	\begin{equation}\label{eq tezine}
				D_{i,j}^\gamma= SI_{i,j}^\gamma+ \frac{h^2}{4}(SI_{i+1,j+1}^\gamma+SI_{i-1,j-1}^\gamma+SI_{i-1,j+1}^\gamma+SI_{i+1,j-1}^\gamma) .
	\end{equation}
	We should point out that although equation \eqref{eq tezine} is similar to equation \eqref{eq 2DsmoothnessFinal}, closest neighbouring 1D indicators are in different positions. 
	Furthermore,
	\begin{equation}\label{eq IzrazOmega2dRavno}				
		\omega_{i,j}^\gamma=\frac{\alpha_{i,j}^\gamma}{\sum_{k=0}^3\alpha_{i,j}^{k\pi/2}},
	\end{equation}
	
	\begin{equation}\label{eq 2dKombinacijaRavno}
		u_{i,j}=\sum_{k=0}^3 \omega_{i,j}^{k\pi/2} p_{i,j}^{k\pi/2}.
	\end{equation}
		Here, $\omega_{i,j}^\gamma$ are non-linear weights dependent on 1D smoothness indicators of nearby polynomials 
and $p_{i,j}^\gamma$ are associated values of interpolation polynomials on stencil $S_{i,j}^{\gamma}$ at the intersecting point $x_{i,j}$.
		The algorithm \ref{alg2} describes the second phase.

	\begin{algorithm}
	\SetAlgoLined
	\footnotesize
	\For {$k=0, \, k< 3,\, k=k+1 $}{
		$\gamma_k=k\pi/2$\\
		\For {$i=1, \,  i<2n-1, \,  i=i+2$}{
			\For {$j=0, \,  j<2m-1, \,  j=j+2$}{
				Calculate $p_{i,j}^{\gamma_k}$  Equation (\ref{eq:p})\\
				Calculate $SI_{i,j}^{\gamma_k}$ Equation (\ref{eq 2Dsmoothness})
			}
		}
		\For {$i=0, \,  i<2n-1, \,  i=i+2$}{
			\For {$j=1, \,  j<2m-1, \,  j=j+2$}{
				Calculate $p_{i,j}^{\gamma_k}$  Equation (\ref{eq:p})\\
				Calculate $SI_{i,j}^{\gamma_k}$, Equation (\ref{eq 2Dsmoothness})
			}
		}
		\For {$i=1, \,  i<2n-1, \,  i=i+2$}{
			\For {$j=0, \,  j<2m-1, \,  j=j+2$}{
				Calculate $\alpha_{i,j}^{\gamma_k}$ Equation (\ref{eq:alpha})\\
				$u_{i,j}=u_{i,j}+\alpha_{i,j}^{\gamma_k} p_{i,j}^{\gamma_k}$
			}
		}
		\For {$i=0, \,  i<2n-1, \,  i=i+2$}{
			\For {$j=1, \,  j<2m-1, \,  j=j+2$}{
				Calculate $\alpha_{i,j}^{\gamma_k}$ Equation (\ref{eq:alpha})\\
				$u_{i,j}=u_{i,j}+\alpha_{i,j}^{\gamma_k} p_{i,j}^{\gamma_k}$
			}
		}			
	}		
	\For {$i=1, \,  i<2n-1, \,  i=i+2$}{
		\For {$j=0, \,  j<2m-1, \,  j=j+2$}{
			$u_{i,j}=\frac{u_{i,j}}{\sum_{k=0}^3 \alpha_{i,j}^{\gamma_k}} $
		}
	}
	\For {$i=0, \,  i<2n-1, \,  i=i+2$}{
		\For {$j=1, \,  j<2m-1, \,  j=j+2$}{
			$u_{i,j}=\frac{u_{i,j}}{\sum_{k=0}^3 \alpha_{i,j}^{\gamma_k}} $
		}
	}
	\caption{{WD} WENO phase 2, $+$ directions }
	\label{alg2}
\end{algorithm}

For polynomial approximation (\ref{eq: Ideal2DPolynomial}) we can state the following proposition.

\begin{prop} \label{proppolinterpol}
Let $v\in C^4(\mathbb{R}^2)$, then the polynomial approximation $q$ from \eqref{eq: Ideal2DPolynomial} based on uniform sampling of $v$ on a rectangular $n\times m$ grid,  satisfy:

\begin{equation}
	v(x_i, y_j)-q(x_i, y_j)=\mathcal{O}(h^4)
\label{eq: Ideal2DPolynomialDif}
\end{equation}

\end{prop}
\begin{proof}
	The proof follows directly from the Taylor expansion of $v$ at $(x_i, y_j)$.
$$
\begin{aligned}
		v(x_i, y_j)-q(x_{i,j})&=\delta\left(v(x_i, y_j)-\sum_{k=0}^1 C_{i,j}^{\gamma_{2k}} p_{i,j}^{\gamma_{2k}}\right)+(1-\delta)\left(v(x_i, y_j)-\sum_{k=0}^1 C_{i,j}^{\gamma_{2k+1}} p_{i,j}^{\gamma_{2k+1}}\right)\\
		&=\mathcal{O}(h^4).
\end{aligned}
$$
\end{proof}

Finally, similar statement can be done for our approximation.
	
\begin{theorem}
Let $v\in C^4(\mathbb{R}^2)$ be sampled on uniform rectangular $n\times m$ grid with constant grid spacing $h$.
Then two phase approximation $u$ defined in \eqref{eq:Vektor}--\eqref{eq 2dKombinacijaRavno}, 
for $(x,y)\in\{ (x_0+ih/2,y_0+jh/2)|0\leq i< 2n, 0\leq j< 2m   \}$ satisfies
$$u(x,y)-v(x,y)=\mathcal{O}(h^4).$$
If $v$ is not smooth everywhere but is at least $C^3$ in some subset of the domain and if at least in one direction the stencil lies in a smooth region, then with $\beta\geq \frac{3}{2}$ and for all $(x,y)$ in this smooth region:
$$u(x,y)-v(x,y)=\mathcal{O}(h^3).$$
\label{th1}
\end{theorem}	
\begin{proof}
At each point of interpolation with index $(i,j)$ there are four second degree polynomials which in two pairs are part of two perpendicular lines containing $(i,j)$, see figures \ref{fig:+} and  \ref{fig:x}. 
The equation \eqref{eq: Ideal2DPolynomialDif} holds for horizontal--vertical directions as well as for slanted directions which can be verified by Taylor expansion of $v$  at $x_{i,j}$. Using the usual WENO argument we should prove that this order of accuracy still holds  when we replace the ideal weights with non-linear WENO weights in smooth regions {\sl i.e.} 
{\setlength\arraycolsep{2pt}
\begin{eqnarray*}
\omega_{i,j}^{\gamma_{2k}}-\delta C_{i,j}^{\gamma_{2k}} & = & \mathcal{O}(h^2), \quad k=0,1, \\
\omega_{i,j}^{\gamma_{2k+1}}-(1-\delta) C_{i,j}^{\gamma_{2k+1}} & = & \mathcal{O}(h^2), \quad k=0,1.
\end{eqnarray*}}%

From the \cite{Liu2009,ARANDIGA2012} and equation \eqref{eq: SIporperties} we can use the properties of 1D WENO smoothness indicators, if $v$ is smooth locally in direction of $\gamma$ :
\begin{equation}
\begin{aligned}\label{eq: SmoothnesOrder}
	SI^{\gamma+\pi}_{i, j}-SI^\gamma_{i, j}= \mathcal{O}\left(h^4\right)\\
	SI^\gamma_{i, j}=\mathcal{O}\left(h^2\right)
\end{aligned}
\end{equation}

\begin{equation}\label{eq: DDiffOrder}
\begin{aligned}
	D^\gamma_{i, j}&=\mathcal{O}\left(h^2\right)\\
	D^{\gamma+\pi}_{i, j}-D^\gamma_{i, j}&= SI^{\gamma+\pi}_{i, j}-SI^\gamma_{i, j}\\
&+\frac{h^2}{8}(SI_{i+1,j+1}^{\gamma+\pi}+SI_{i-1,j-1}^{\gamma+\pi}+SI_{i-1,j+1}^{\gamma+\pi}+SI_{i+1,j-1}^{\gamma+\pi})\\
&-\frac{h^2}{8}(SI_{i+1,j+1}^\gamma+SI_{i-1,j-1}^\gamma+SI_{i-1,j+1}^\gamma+SI_{i+1,j-1}^\gamma)\\
&=\mathcal{O}\left(h^4\right)+ \frac{h^2}{4}\mathcal{O}\left(h^4\right)
\end{aligned}
\end{equation}

\begin{equation}\label{eq: D}
	\begin{aligned}
		\frac{\frac{1}{\left(\varepsilon_h+D^\gamma_{i, j}\right)^\beta}-\frac{1}{\left(\varepsilon_h+D^{\gamma+\pi}_{i, j}\right)^\beta}}{\frac{1}{\left(\varepsilon_h+D^{\gamma+\pi}_{i, j}\right)^\beta}} & =\left(\frac{\varepsilon_h+D^{\gamma+\pi}_{i, j}}{\varepsilon_h+D^\gamma_{i, j}}\right)^\beta-1 \\
		& =\left(\frac{\varepsilon_h+D^{\gamma+\pi}_{i, j}}{\varepsilon_h+D^\gamma_{i, j}}-1\right) \sum_{l=0}^\infty {{\beta}\choose{l+1}} \left(\frac{\varepsilon_h+D^{\gamma+\pi}_{i, j}}{\varepsilon_h+D^\gamma_{i, j}}-1\right)^l \\
		& =\frac{D^{\gamma+\pi}_{i, j}-D^\gamma_{i, j}}{K_\varepsilon h^2+\mathcal{O}\left(h^2\right)} \sum_{l=0}^\infty {{\beta}\choose{l+1}} \left(\frac{D^{\gamma+\pi}_{i, j}-D^\gamma_{i, j}}{K_\varepsilon h^2+\mathcal{O}\left(h^2\right)}\right)^l \\
		& =\frac{\mathcal{O}\left(h^{4}\right)}{K_\varepsilon h^2+\mathcal{O}\left(h^2\right)} \sum_{l=0}^{\infty} {{\beta}\choose{l+1}}\left(\frac{\mathcal{O}\left(h^{4}\right)}{K_\varepsilon h^2+\mathcal{O}\left(h^2\right)}\right)^l\\
		& =\mathcal{O}\left(h^2\right) \mathcal{O}\left(1\right) ,
	\end{aligned}
\end{equation}
where the series in equation \eqref{eq: D} converges near zero.
Therefore

\begin{equation}
	\begin{aligned}
	\frac{1}{\left(\varepsilon_h+D^\gamma_{i, j}\right)^\beta}&=\frac{1}{\left(\varepsilon_h+D^{\gamma+\pi}_{i, j}\right)^\beta}\left(1+\mathcal{O}\left(h^2\right)\right), \quad \forall(i, j),\gamma \\	
	\alpha^{\gamma+\pi}_{i, j}&=\frac{0.5\left(1+\mathcal{O}\left(h^2\right)\right)}{\left(\varepsilon_h+D^{\gamma}_{i, j}\right)^\beta}, \\
	\alpha^\gamma_{i, j}+\alpha^{\gamma+\pi}_{i, j}&=\frac{\left(1+\mathcal{O}\left(h^2\right)\right)}{\left(\varepsilon_h+D^{\gamma}_{i, j}\right)^\beta},\\	
	\delta&=\frac{\alpha^\gamma_{i, j}+\alpha^{\gamma+\pi}_{i, j}}{\alpha^\gamma_{i, j}+\alpha^{\gamma+\pi}_{i, j}+\alpha^{\gamma+\pi/4}_{i, j}+\alpha^{\gamma+3\pi/4}_{i, j}}\\
	\frac{\omega_{i,j}^\gamma}{\delta}=&\frac{\alpha^\gamma_{i, j}}{\alpha^\gamma_{i, j}+\alpha^{\gamma+\pi}_{i, j}} =\frac{0.5 /\left(\varepsilon_h+D^{\gamma}_{i, j}\right)^\beta}{\left(1+\mathcal{O}\left(h^2\right)\right) /\left(\varepsilon_h+D^{\gamma}_{i, j}\right)^\beta}
	=\frac{0.5}{1+\mathcal{O}\left(h^2\right)} \\
	& = 0.5\left(1+\mathcal{O}\left(h^2\right)\right).
	\end{aligned} \label{defdelta}
\end{equation}		

Similarly, we get
$$
\begin{aligned}
\frac{\omega_{i,j}^\gamma}{\delta}-0.5&=\mathcal{O}\left(h^2\right),\:
\frac{\omega_{i,j}^{\gamma+\pi}}{\delta}-0.5=\mathcal{O}\left(h^2\right)\\
\frac{\omega_{i,j}^{\gamma+\pi/4}}{1-\delta}-0.5&=\mathcal{O}\left(h^2\right),\:
\frac{\omega_{i,j}^{\gamma+3\pi/4}}{1-\delta}-0.5=\mathcal{O}\left(h^2\right)
\end{aligned}
$$

If, $v$ has a singularity inside stencil of $D^\gamma$, then $D^{\gamma}_{i, j} =\mathcal{O}\left(1\right)$, whereas $D^{\gamma}_{i, j}=O\left(h^2\right)$ otherwise, then
$$
\alpha^{\gamma}_{i, j}= \begin{cases}\mathcal{O}(1), & v \text { not smooth }, \\ \mathcal{O}\left(h^{-2 \beta}\right), & v \text { smooth } ,\end{cases}
$$
therefore $\sum_{k=0}^{3} \alpha^{\gamma_k}_{i,j}=\mathcal{O}\left(h^{-2 \beta}\right)$ and $\omega_{i,j}^{\gamma}=\mathcal{O}\left(h^{2 \beta}\right)$ if $v$ is not smooth in direction $\gamma$ but smooth in at least one of the possible directions. If we denote $\mathcal{S}=\{k| \textrm{v is smooth in}$ $\textrm{direction }\gamma_k$\} then
$$
\begin{aligned}
	v_{i,j}-u_{i,j} & =\sum_{k\notin \mathcal{S}} \omega_{i,j}^{\gamma_k}\left(v_{i,j}-p^{\gamma_k}_{i,j}\right) +\sum_{k\in \mathcal{S}} \omega_{i,j}^{\gamma_k}\left(v_{i,j}-p^{\gamma_k}_{i,j}\right)\\
	& =\sum_{k \notin \mathcal{S}} \mathcal{O}\left(h^{2 \beta}\right) \mathcal{O}(1)+\sum_{k \in \mathcal{S}} \mathcal{O}(1) \mathcal{O}\left(h^3\right)=\mathcal{O}\left(h^{\min (3, 2 \beta)}\right) .
\end{aligned}
$$

If $v$ is smooth locally in all directions, {WD} WENO approximation almost reduces to  a weighted average of interpolating cubic polynomials. We can write the proof for directions $\gamma_k=\gamma_0+k\pi/2$ and split the indices in two sets   $\mathcal{S}_1=\{0,2\}$ and  $\mathcal{S}_2=\{1,3\}$, then from Proposition \ref{proppolinterpol} we have
$$
\begin{aligned}
	v_{i,j}-u_{i,j} & =	v_{i,j}-q_{i,j}+q_{i,j}-u_{i,j}=v_{i,j}-\delta\sum_{k\in \mathcal{S}_1} C_{i,j}^{\gamma_k}p^{\gamma_k}_{i,j}-(1-\delta)\sum_{k\in \mathcal{S}_2} C_{i,j}^{\gamma_k}p^{\gamma_k}_{i,j} \\
	&+\delta\sum_{k\in \mathcal{S}_1} C^{\gamma_k}_{i,j}\left(p_{i,j}^{\gamma_k}-v_{i,j}\right)+(1-\delta)\sum_{k\in \mathcal{S}_2} C^{\gamma_k}_{i,j}\left(p_{i,j}^{\gamma_k}-v_{i,j}\right) \\
	&+\sum_{k\in \mathcal{S}} \omega^{\gamma_k}_{i,j}\left(v_{i,j}-p_{i,j}^{\gamma_k}\right)\\
	& =	\mathcal{O}\left(h^4\right)- \delta\sum_{k\in \mathcal{S}_1} \left(p^{\gamma_k}_{i,j}-v_{i,j}\right)\left(\omega_{i,j}^{\gamma_k}/\delta-0.5\right)\\
	&-(1-\delta)\sum_{k\in \mathcal{S}_2} \left(p^{\gamma_k}_{i,j}-v_{i,j}\right)\left(\omega_{i,j}^{\gamma_k}/(1-\delta)-0.5\right)\\
	& = \mathcal{O}\left(h^4\right)+\delta \mathcal{O}(h^3) \mathcal{O}\left(h^2\right)+(1-\delta) \mathcal{O}(h^3) \mathcal{O}\left(h^2\right)=\mathcal{O}\left(h^4\right) .
\end{aligned}
$$
Finally, if some of the values used in the interpolation are calculated in previous phase in the smooth regions, the error is small enough to not interfere with the order of interpolation.
\vspace{-1em}
\flushright{Q.E.D.}
\end{proof}	

\subsection{Arbitrary resolution interpolation} \label{sec: Phase3}
There is a simple way of obtaining an interpolation method which can be applied to get 2D interpolation with an arbitrary resolution, relaying on the 1D interpolation from  section \ref{sec Interpolation} which was used in \cite{Mulet_2010}. 
Based on local WENO interpolation function \eqref{eq TezinskiPolinom1d} it useful to introduce 1D WENO interpolation spline which interpolates function values $V= (v_0,\ldots v_{n-1})$ on uniformly spaced grid $0,\ldots,n-1$ which we want to evaluate on non-integer valued grid $0,\ldots, x_{n^*-1}$ which satisfies $x_{n^*-1}=n-1$:
$$ \label{eq wenoSpline}
w(V,x)=u_i((v_{i-1},v_i,v_{i+1},v_{i+2}),x), \quad x\in[x_i,x_{i+1}]. 
$$ 
The images are usually represented as matrices $A:=(a_{i,j})$, where $a_{i,j}$ is located in the $i$-th row and $j$-th column. Although our notation of $x_{i,j}$ is little bit different, it actually doesn't matter, so for the simplicity we will denote with $a_{i,j}$ the original value of the pixel at the position $x_{i,j}$, and $A\in\mathbb{R}^{n\times m}$. Now we can construct a tensor product interpolation function for given matrix:
\begin{equation}
\begin{aligned}	
	z_k&=w(A(:,k),y), \quad k=0,\ldots,m-1,\\
	Z&=(z_0,\dots,z_{m-1}),\\
	U(x,y)&=w(Z,x).
	\label{eq: tensor}
\end{aligned}	
\end{equation}

\begin{figure}
	\begin{center}
		{
			\subfigure[Input image]{\includegraphics[width=0.3\textwidth]{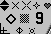}}
			\subfigure[Tensor WENO ]{\includegraphics[width=0.3\textwidth]{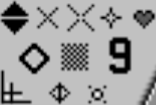}}
			\subfigure[{WD} WENO ]{\includegraphics[width=0.3\textwidth]{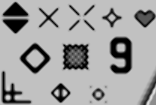}}
		}			
		\caption{Pixel art image interpolation with scale facotor 3.}
		\label{fig: scale3}
	\end{center}
\end{figure}

This tensor WENO algorithm provides the advantage of separable interpolation, allowing the interpolation process to be broken down into one-dimensional interpolations along each dimension. 
Its complexity is the same as {WD} WENO interpolation and satisfies both the non-oscillatory properties and the accuracy requirements outlined in Theorem \ref{th1} for {WD} WENO. 
However, the algorithm falls short in achieving the desired outcomes for images, as evident in Figure \ref{fig: scale3} b), where the image was upsampled by a factor of $3$. 
Notably, the tensor WENO algorithm introduces a staircase effect on slanted lines. 
On a positive note, it's worth mentioning that this algorithm is adaptable to any resolution. 
Furthermore, the method produces minimal artifacts during downsampling, making it a viable component of {WD} WENO algorithm (not shown).

If the resolution of the interpolated image precisely matches
\begin{equation}
	(2^k(n-1)+1)\times(2^k(m-1)+1)
	\label{eq:zoomfactor}
\end{equation}
for some integer $k$, then this image can be interpolated by successively applying the {WD} WENO algorithm with two phases exactly $k$ times. 
In cases where the resolution does not satisfy equation \eqref{eq:zoomfactor}, we initially interpolate an image with a larger magnification that satisfies the equation \eqref{eq:zoomfactor}. 
The resulting image is then downsampled to the exact resolution using a tensor WENO algorithm  \eqref{eq: tensor}. 
The outcomes of this algorithm can be observed in Figure \ref{fig: scale3} c), where the image was initially upsampled to dimensions $(4n-3)\times(4m-3)$ and subsequently downsampled to dimensions $3n\times3m$ using tensor WENO.

The images generated through the aforementioned {WD} WENO algorithm can be employed for upsampling or downsampling to arbitrary resolutions, preserving all the  effects originally designed for image resolution "doubling". 
It is evident that "doubling" the image resolution can be executed faster and with a reduced memory load compared to a more general scaling factor. 
Consequently, the numerical tests will be concentrated on resolutions that fulfills the condition specified in \eqref{eq:zoomfactor}.

\section{Numerical tests} 
	\label{sec:numtests}
	
To evaluate the behavior of the {WD} WENO interpolation, we will use selected images to illustrate and demonstrate some important properties of the proposed algorithm and compare it with other similar methods commonly used for upsampling images.
The proposed method will be tested on images with sharp color transitions and on standard real-world test images taken with a digital camera.

The quality of the upscaled images is measured using the Peak Signal-to-Noise Ratio (PSNR) \cite{Gonzales} and the Mean Structural Similarity Index (MSSIM) \cite{Wang}, which are commonly used to measure image quality.
In addition to these measurements, images must be visually inspected to rule out the possibility of local numerical artifacts that do not affect PSNR and MSSIM.

Unless stated otherwise, the free parameter $\beta$ used in Eq. \eqref{eq:alpha} is set to $\beta = 2$. Furthermore, the upsampled images will satisfy equation \eqref{eq:zoomfactor} but for simplicity we will designate as scale factors $d=2,4,8,16$.
	
\subsection{Convergence test}
We verify numerically that the order of the error for the interpolation proposed in this paper is $\mathcal{O}(h^4)$ for sufficiently smooth functions. In this experiment, we interpolate the function
$$f(x,y)=\frac{1}{x^2+y^2+1}$$
on a square domain $[-1,1] \times [-1,1]$ using a regular rectangular grid with grid spacing $h=2^i, i=0,\ldots10$. Table \ref{tab:convergence} shows $L_2$ and $L_{inf}$ errors as well as the numerical convergence rate of the {WD} WENO algorithm with exponent $\beta=1$. The convergence rates in both norms are very close to the predicted rates in the Theorem \ref{th1}
	
	 \begin{table}
		\centering\footnotesize
		\begin{tabular}{ccccc}
			\hline
			$h$ & $L_{inf}$ error &  $\mathcal{O}_{inf} (h^*)$ &  $L_2$ error  & $\mathcal{O}_2 (h^*)$ \\
			\hline
1.00E+00	&	6.32E-01	&		&	3.24E-01	&		\\
5.00E-01	&	1.33E-01	&	2.25	&	7.20E-02	&	2.17	\\
2.50E-01	&	4.29E-02	&	1.63	&	1.50E-02	&	2.26	\\
1.25E-01	&	6.51E-03	&	2.72	&	1.58E-03	&	3.25	\\
6.25E-02	&	5.15E-04	&	3.66	&	1.10E-04	&	3.85	\\
3.13E-02	&	3.47E-05	&	3.89	&	6.77E-06	&	4.02	\\
1.56E-02	&	2.36E-06	&	3.88	&	4.14E-07	&	4.03	\\
7.81E-03	&	1.56E-07	&	3.92	&	2.56E-08	&	4.02	\\
3.91E-03	&	1.01E-08	&	3.95	&	1.59E-09	&	4.01	\\
1.95E-03	&	6.42E-10	&	3.97	&	9.95E-11	&	4	\\
9.77E-04	&	4.05E-11	&	3.99	&	6.22E-12	&	4	\\
			\hline			
		\end{tabular}
		\caption{Rate of convergence approximated from $L_{inf}$ and $L_2$ errors for a sequence of embedded numerical meshes $\mathcal{O}(h^*)=\log_2 \frac{Err_i}{Err_{i+1}}$}
		\label{tab:convergence}
	\end{table}
	\subsection{Discontinuity test}
\begin{figure}
	\centering
	\subfigure[Discontinuous function]{\includegraphics[width=0.4\textwidth]{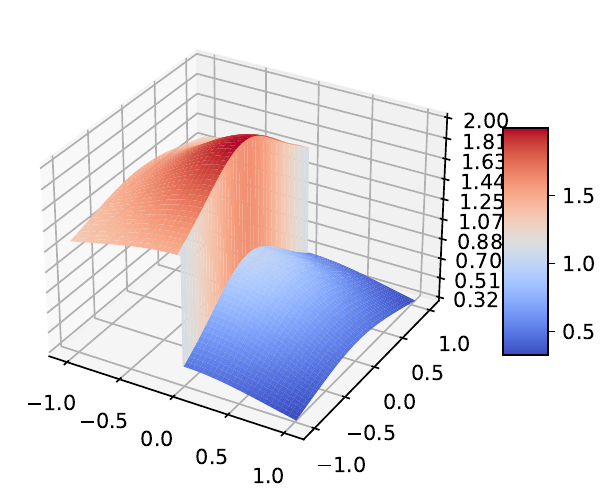}}
	\subfigure[Distribution of relative errors ]{\includegraphics[width=0.4\textwidth]{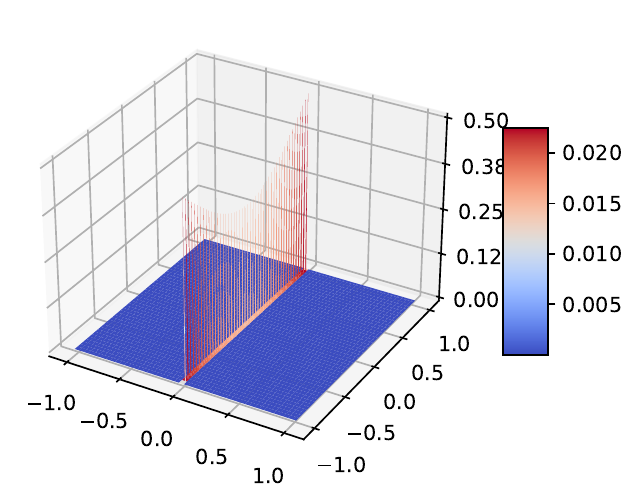}}
	\caption{Interpolated discontinuous function using {WD} WENO $\beta=2$}
	\label{figSurface}
	\end{figure}

 In this experiment we interpolate the function
 \begin{equation}
 	f(x,y)=
 	\begin{cases}
 		\frac{1}{x^2+y^2+1}+1,& \text{if } x < 0\\
 		\frac{1}{x^2+y^2+1}, & \text{otherwise}
 	\end{cases}
 \end{equation}
 on a square domain $[-1,1] \times [-1,1]$ using a regular rectangular grid with grid spacing $h=2^i, i=0,\ldots10$. Although we can not expect convergence in $L_{inf}$, we can observe the behavior of the interpolation algorithm and the relative error distribution in Figure \ref{figSurface}. The errors are significantly larger on the left-hand side of the discontinuity and the convergence in $L_{2}$ is linear at best, which can be seen in Table \ref{tab:convergenceDiscontinuity}, since close to the discontinuity, there are some points where the function $f$ is not smooth in all four directions at some phase. Therefore, no condition of Theorem \ref{th1} is fulfilled. If we measure the errors only on the right-hand side of the discontinuity ($x\in[0,1]$), we can expect convergence of order three according to theorem \ref{th1}, since only a part of the stencil of {WD} WENO method is in the smooth region. In Table \ref{tab:convergenceDiscontinuity}, the measured convergence in the last two columns is even better than expected, confirming Theorem \ref{th1}. For the parameter $\beta=1$, the order is less than three, which is expected according to theorem \ref{th1} since $\beta <1.5$. Moreover, for the parameter $\beta=0$, the order of the approximation measured in the right-hand part of the domain decreases to first order, since the approximation in this case can not avoid the discontinuity. This demonstrates the Theorem \ref{th1} and the adaptive behavior of the presented {WD} WENO approximation.
 
	 \begin{table}
	\centering\footnotesize
	\begin{tabular}{c|cc|cc|cc|cc}
		 &\multicolumn{2}{c}{$\beta=2,x\in[-1,1]$}& \multicolumn{2}{c}{$\beta=0,x\in[0,1]$}  &\multicolumn{2}{c}{$\beta=1,x\in[0,1]$}&\multicolumn{2}{c}{$\beta=2,x\in[0,1]$}  \\
	\hline	
	h	&	$\mathcal{O}_{inf} (h^*)$	&	$\mathcal{O}_2 (h^*)$	&	$\mathcal{O}_{inf} (h^*)$	&	$\mathcal{O}_2 (h^*)$	&	$\mathcal{O}_{inf} (h^*)$	&	$\mathcal{O}_2 (h^*)$&	$\mathcal{O}_{inf} (h^*)$	&	$\mathcal{O}_2 (h^*)$\\
	\hline
1.00E+00	&		&		&		&		&		&		&		&		\\
5.00E-01	&	-1.58	&	-0.66	&	0.8	&	1.77	&	0.0	&	1.0	&	-0.05	&	0.72	\\
2.50E-01	&	-0.01	&	0.36	&	0.45	&	0.57	&	1.2	&	1.6	&	1.6	&	2.02	\\
1.30E-01	&	0.14	&	0.55	&	0.06	&	0.56	&	2.2	&	2.6	&	2.8	&	3.19	\\
6.30E-02	&	0.1	&	0.54	&	0.01	&	0.51	&	1.9	&	2.9	&	3.52	&	3.95	\\
3.10E-02	&	0.03	&	0.51	&	0	&	0.5	&	2.0	&	2.6	&	3.77	&	4.19	\\
1.60E-02	&	0.01	&	0.5	&	0	&	0.5	&	2.0	&	2.5	&	3.55	&	4.14	\\
7.80E-03	&	0	&	0.5	&	0	&	0.5	&	2.0	&	2.5	&	3.73	&	4.06	\\
3.90E-03	&	0	&	0.5	&	0	&	0.5	&	2.0	&	2.5	&	3.9	&	4.03	\\
	\hline			
	\end{tabular}
	\caption{Rate of convergence of {WD} WENO algorithm measured on different parts of the domain}
	\label{tab:convergenceDiscontinuity}
	\end{table}
	
\subsection{Influence of the free parameter $\beta$}
We have seen in the previous section the influence of the parameter $\beta$ on the convergence of the interpolation algorithm for smooth and piece-wise smooth functions.
In this test, we want to get a qualitative evaluation and direct visual confirmation of how the quality of the upsampled image depends on the choice of the free parameter $\beta$ \eqref{eq:alpha}.
The input image \ref{fig:head} was upsampled using our method, as described in the Subsection \ref{sec:natural}.

\begin{figure}
	\begin{center}
		\subfigure[Input image]{\includegraphics[width=0.30\textwidth]{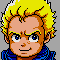}} \quad 
		\subfigure[$\beta=0$]{\includegraphics[width=0.30\textwidth]{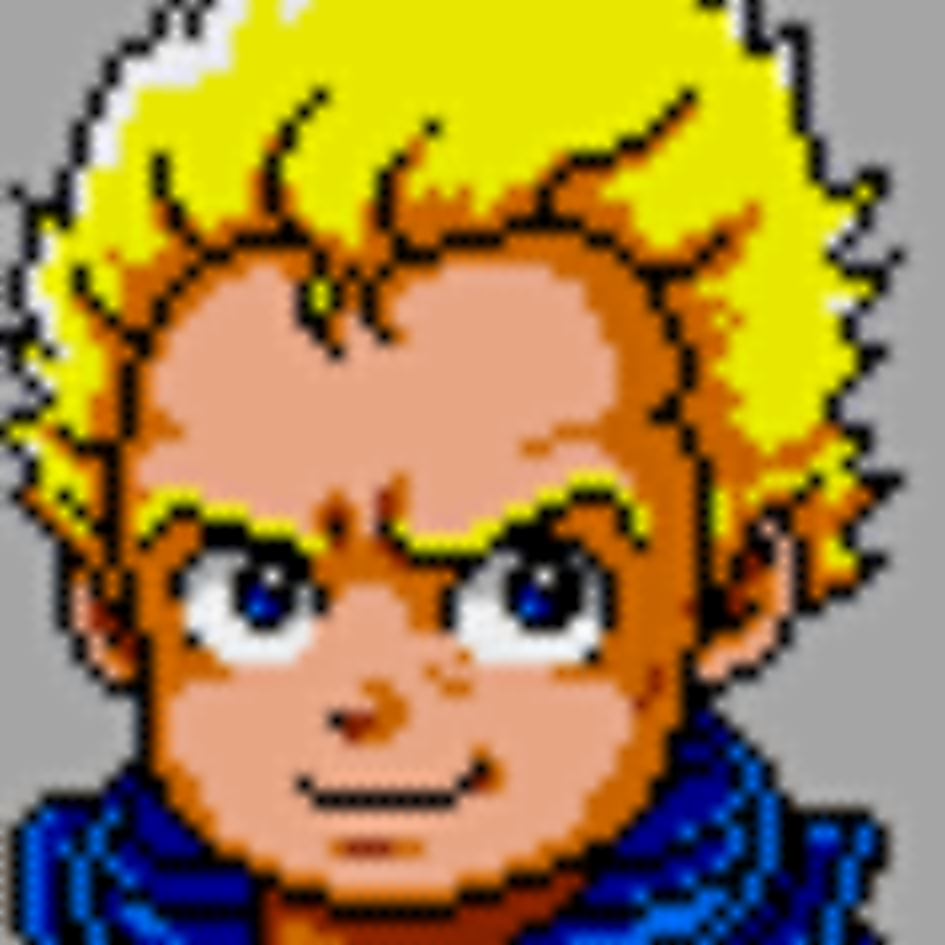}} \quad 
		\subfigure[$\beta=0.1$]{\includegraphics[width=0.30\textwidth]{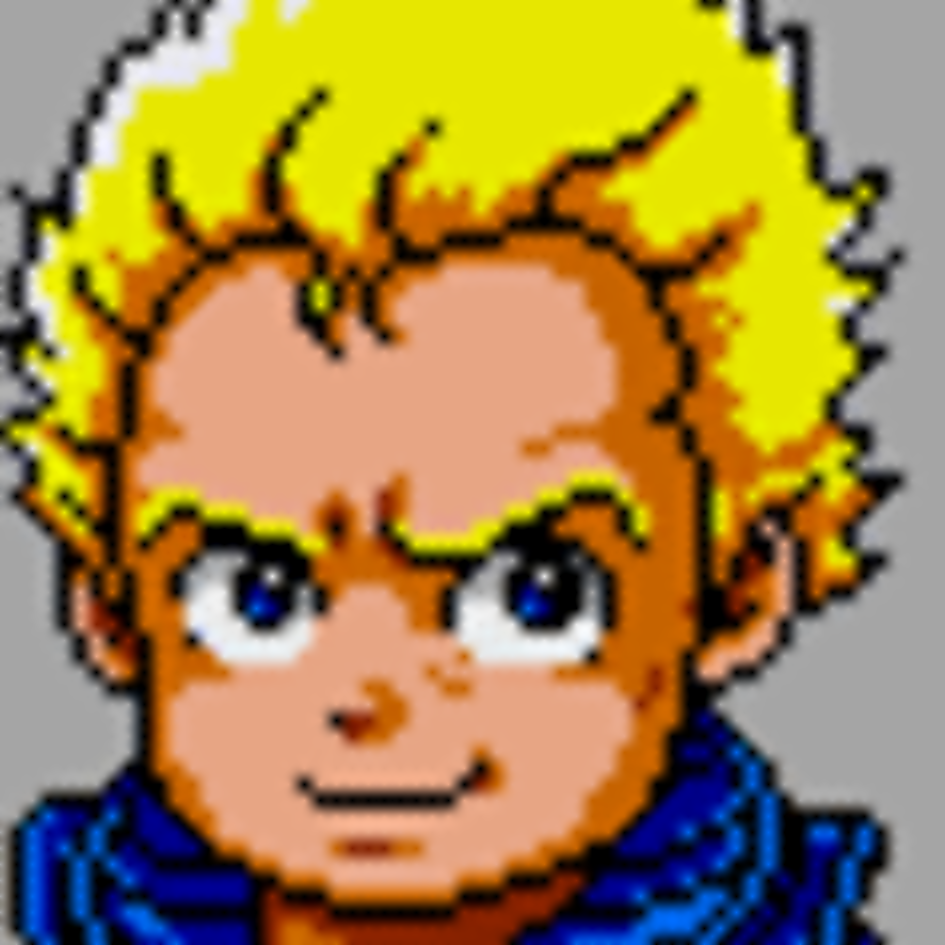}} \quad
		\subfigure[$\beta=1$]{\includegraphics[width=0.30\textwidth]{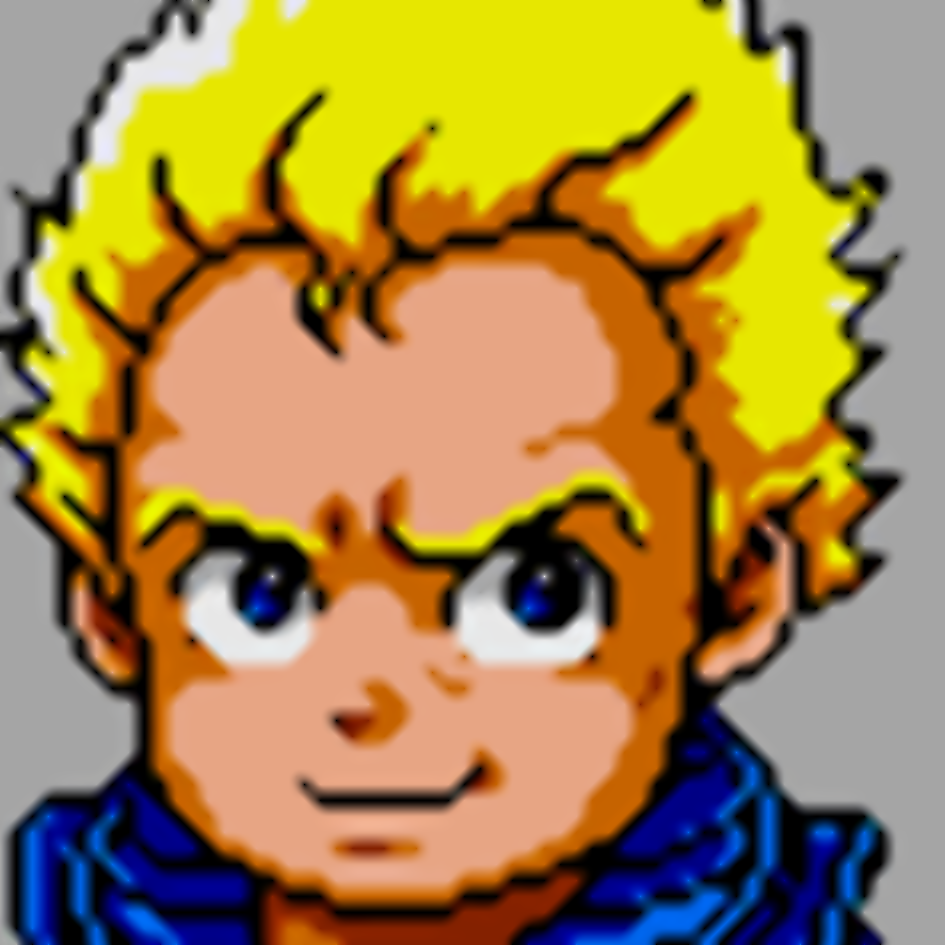}} \quad
		\subfigure[$\beta=2$]{\includegraphics[width=0.30\textwidth]{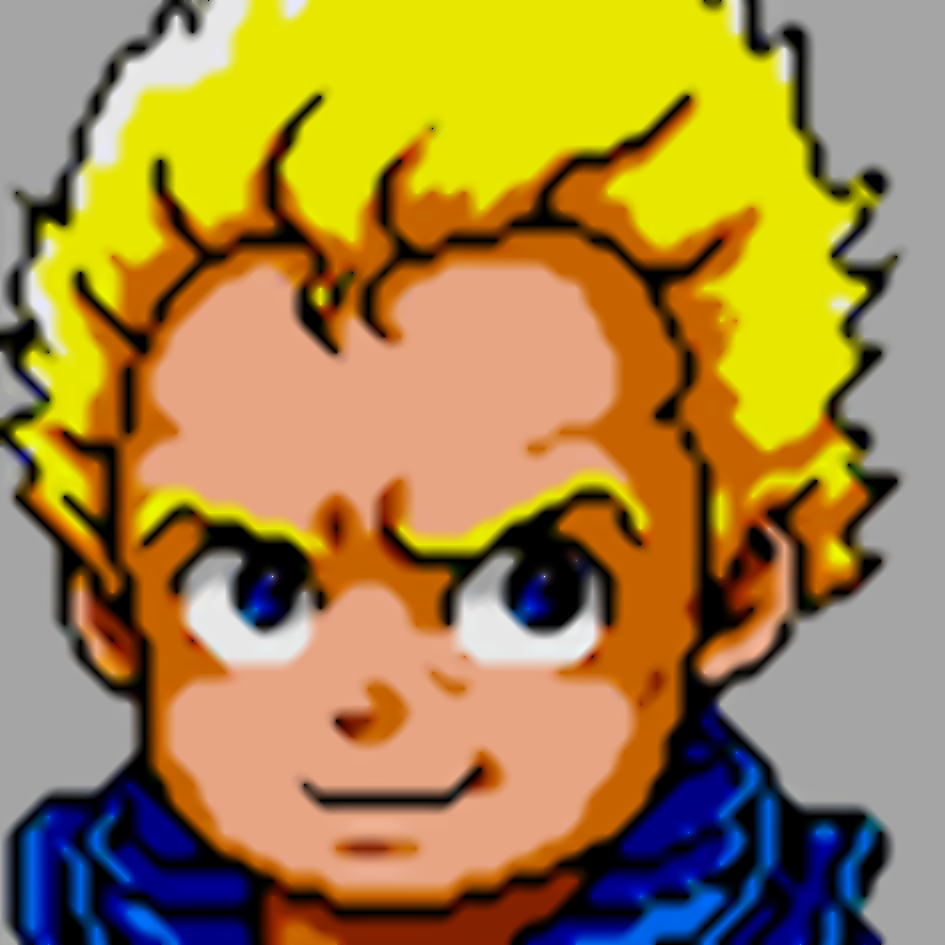}} \quad
		\subfigure[$\beta=3$]{\includegraphics[width=0.30\textwidth]{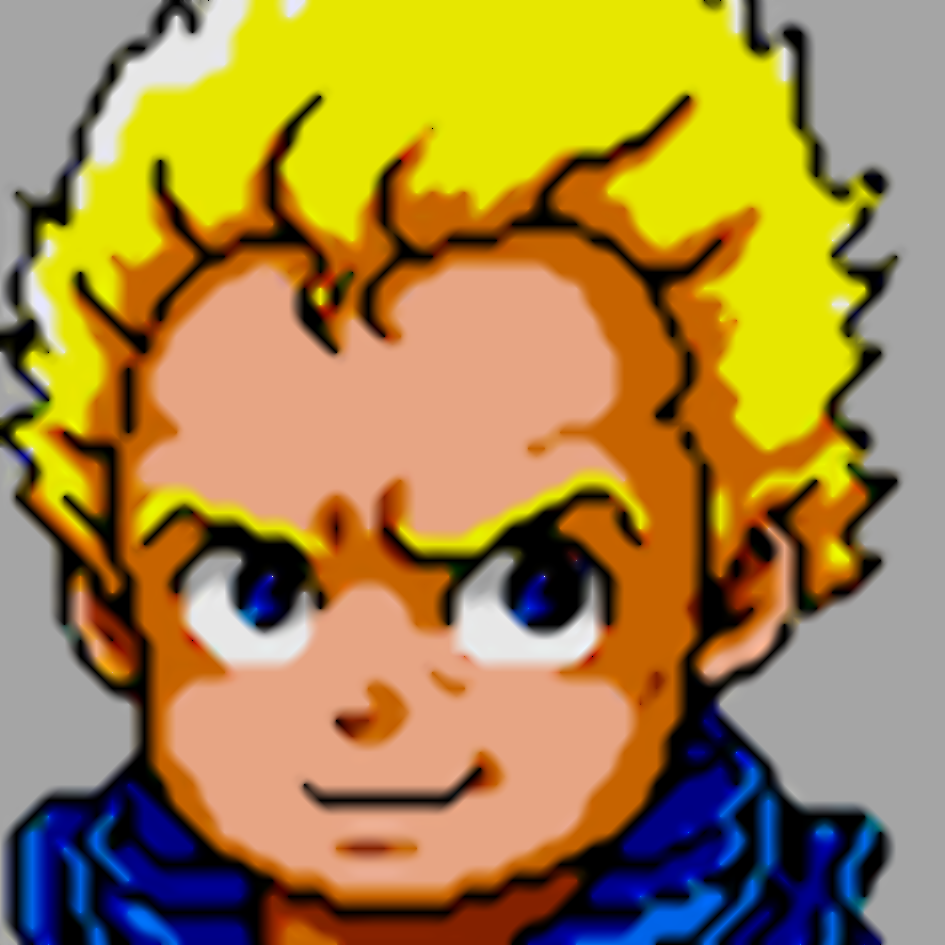}} \quad
		\caption{Influence of $\beta$ exponent on interpolation, scale $d=16$}
		\label{fig:smiley:beta}
	\end{center}
\end{figure}

The resulting MSSIM and PSNR values showed no significant deviations for different $\beta$ values, in contrast to visual inspection of the scaled images.
The effect is clearly seen in Figure \ref{fig:smiley:beta}, where the input image was scaled for $d=16$.
For this reason, the free parameter $\beta=2$ was chosen from this point on.
	
	\subsection{Natural images}\label{sec:natural}
	Our first numerical test was conducted with the standard Kodak  Lossless True Color Image Suite \cite{Kodak}, a set of 24 natural color images.  
	
	The performance of our {WD} WENO algorithm was measured for two scaling factors $d=2, 4$.
	Each image is first cropped by the last row and column to obtain the referent image for $d=2$. 
	Similarly, in order to obtain 
	the referent image for $d=4$, original image was cropped by last three rows and columns.

	Then the images are downsampled according to the expression
	$$\mathrm{downsampled}_{sz} = (\mathrm{referent}_{sz}- 1)/d+1. $$
	The downsampled image was then upsampled with tested methods according to the expression 
	$$\mathrm{scaled}_{sz} = d(\mathrm{downsampled}_{sz} - 1)+1.$$
	
	Since the results are very consistent across all Kodak suit images, we will present the average result rather than picking out a representative image (Table \ref{table:kodims:scale24}).
	
	\begin{table}
		\begin{center}
			\caption{PSNR and MSSIM on the Kodim image suite for d=2, 4}
			\label{table:kodims:scale24}
			\begin{tabular}{lrrrr|rrrr}
				\hline
				& \multicolumn{4}{c}{$d=2$}& \multicolumn{4}{c}{$d=4$} \\
				Method          & PSNR      &Rank   &MSSIM      &Rank   &PSNR       &Rank   &MSSIM      &Rank  \\
				\hline
				{WD} WENO &  27.6841  &  2  &  0.8106  &  2  &  23.4737  &  5  &  0.6248  &  1 \\
				\hline
				Box &  25.6948  &  31  &  0.7519  &  31  &  22.1150  &  31  &  0.5635  &  31 \\
				Triangle &  27.2094  &  22  &  0.7899  &  25  &  23.3565  &  9  &  0.6107  &  9 \\
				\hline
				Catrom &  27.3542  &  13  &  0.8029  &  16  &  23.0967  &  14  &  0.6092  &  11 \\
				Cubic &  26.6882  &  29  &  0.7557  &  29  &  23.5256  &  1  &  0.6073  &  14 \\
				Gaussian &  26.9763  &  25  &  0.774  &  27  &  23.4976  &  4  &  0.6112  &  8 \\
				Hermite &  27.1648  &  23  &  0.7928  &  22  &  23.1746  &  13  &  0.6063  &  17 \\
				Lagrange &  27.3729  &  5  &  0.8015  &  18  &  23.1895  &  12  &  0.6114  &  7 \\
				Mitchell &  27.2609  &  20  &  0.792  &  23  &  23.3668  &  8  &  0.6131  &  4 \\
				Robidoux &  27.2363  &  21  &  0.7901  &  24  &  23.3931  &  7  &  0.6132  &  3 \\
				RobidouxSharp &  27.2945  &  19  &  0.7948  &  21  &  23.3202  &  10  &  0.6127  &  5 \\
				Spline &  26.6882  &  29  &  0.7557  &  29  &  23.5256  &  1  &  0.6073  &  14 \\
				\hline
				Bartlett &  27.3627  &  11  &  0.8052  &  13  &  22.9317  &  21  &  0.6032  &  25 \\
				Blackman &  27.3741  &  3  &  0.8056  &  10  &  22.9798  &  18  &  0.6063  &  17 \\
				Bohman &  27.3734  &  4  &  0.8055  &  12  &  22.9868  &  17  &  0.6066  &  16 \\
				\hline
				Cosine &  27.3541  &  14  &  0.8059  &  6  &  22.8947  &  25  &  0.6031  &  26 \\
				Hamming &  27.3554  &  12  &  0.8058  &  7  &  22.8901  &  26  &  0.6026  &  27 \\
				Hann &  27.3636  &  8  &  0.8057  &  9  &  22.921  &  24  &  0.6036  &  23 \\
				Jinc &  27.1287  &  24  &  0.7819  &  26  &  23.4288  &  6  &  0.6094  &  10 \\
				Kaiser &  27.3682  &  7  &  0.8058  &  7  &  22.9368  &  20  &  0.6045  &  22 \\
				Lanczos &  27.3631  &  9  &  0.806  &  4  &  22.925  &  22  &  0.6047  &  20 \\
				Lanczos2 &  27.352  &  16  &  0.8029  &  16  &  23.0892  &  15  &  0.6088  &  12 \\
				Lanczos2Sharp &  27.3146  &  18  &  0.8039  &  15  &  22.9756  &  19  &  0.6055  &  19 \\
				LanczosRadius &  27.3631  &  9  &  0.806  &  4  &  22.925  &  22  &  0.6047  &  20 \\
				LanczosSharp &  27.3522  &  15  &  0.8065  &  3  &  22.8812  &  27  &  0.6035  &  24 \\
				Parzen &  27.3725  &  6  &  0.8051  &  14  &  23.0122  &  16  &  0.6074  &  13 \\
				Point &  25.6948  &  31  &  0.7519  &  31  &  22.115  &  31  &  0.5635  &  31 \\
				Quadratic &  26.9472  &  28  &  0.7715  &  28  &  23.5174  &  3  &  0.6115  &  6 \\
				Sinc &  26.9607  &  26  &  0.7967  &  19  &  22.3428  &  29  &  0.5791  &  29 \\
				SincFast &  26.9607  &  26  &  0.7967  &  19  &  22.3428  &  29  &  0.5791  &  29 \\
				Welch &  27.347  &  17  &  0.8056  &  10  &  22.8776  &  28  &  0.6021  &  28 \\
				\hline
				GCS &  27.7598  &  1  &  0.8166  &  1  &  23.2188  &  11  &  0.6214  &  2 \\
				\hline
			\end{tabular}
		\end{center}
	\end{table}
	
The methods in the tables are grouped by filter type.
The first method is the proposed {WD} WENO method.
The rest of the tested methods are implemented in ImageMagick software \cite{imagemagick-web-2018}, except for the Geometric Contour Stencils (GCS) \cite{getreuergeomcontstenc} method.

The second group of methods are the cubic interpolation methods, which provide relatively fast execution times. These methods may exhibit ringing effects or introduce numerical diffusion to solve this problem.

The box and triangle filters are simple, fast interpolation methods that produce strong aliasing effects but no ringing effects.

Windowed--sinc filters are a separate group of filters that have higher computational costs compared to the interpolation methods and are considered the best filters for use with real images. In general, windowed--sinc filters produce stronger ringing effects, but some of them, such as the Lanczos Sharp filter, are usually recommended for upsampling images with minimal blur and ringing.

The last group shows the Geometric Contour Stencils method. This is a computationally intensive method that estimates the geometric properties of contours in an image that affect the interpolation process. By incorporating information about contour tangents and curvature, the method aims to produce interpolated images that preserve the structure and sharpness of the contours.

The results in Table \ref{table:kodims:scale24} show that the images generated by the {WD} WENO method achieve very good results in the PSNR and MSSIM measures.

\begin{figure}
		\begin{center}
			\subfigure[Input image]{\includegraphics[width=0.30\textwidth]{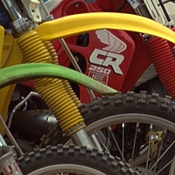}} \quad 
			\subfigure[Catrom]{\includegraphics[width=0.30\textwidth]{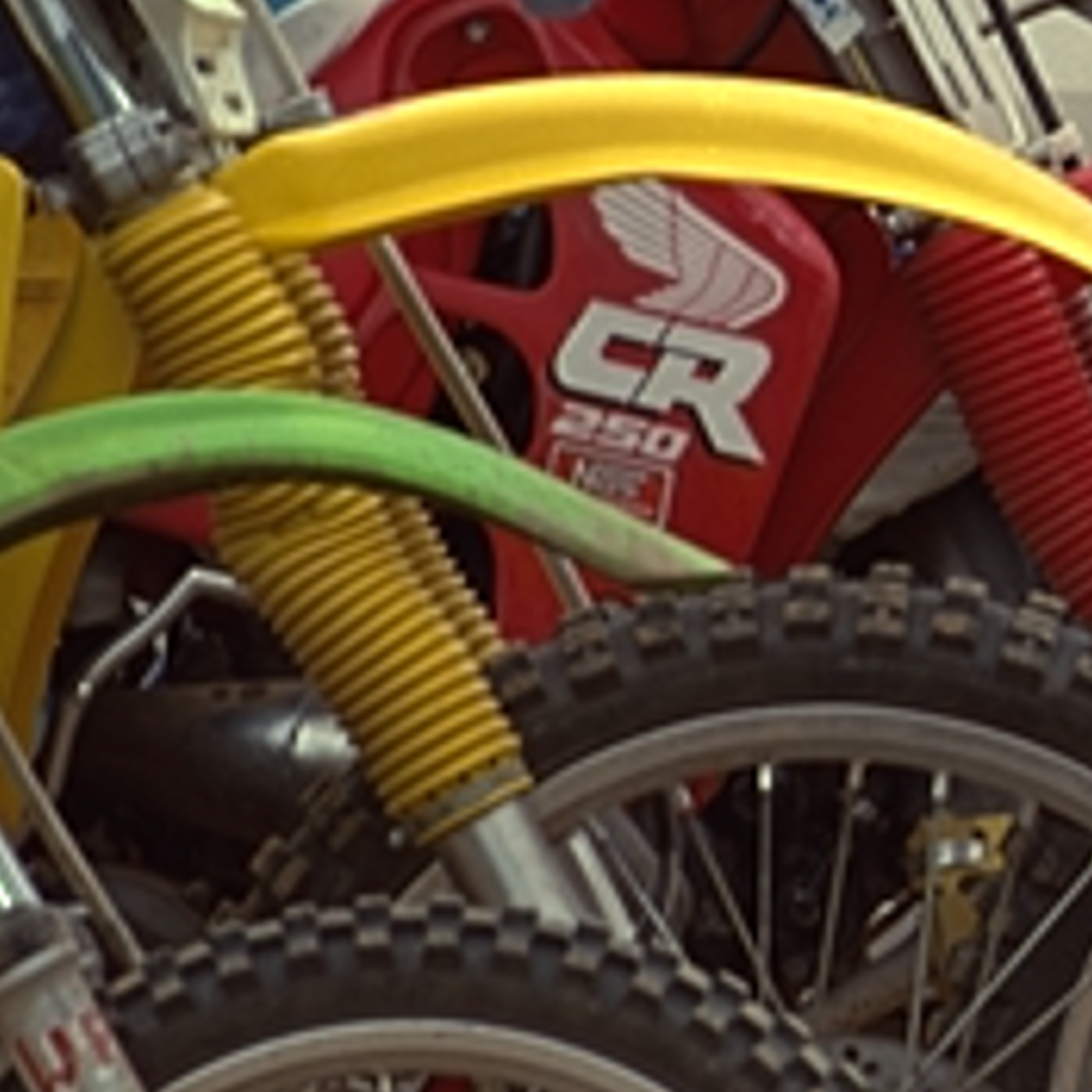}} \quad 
			\subfigure[Robidoux]{\includegraphics[width=0.30\textwidth]{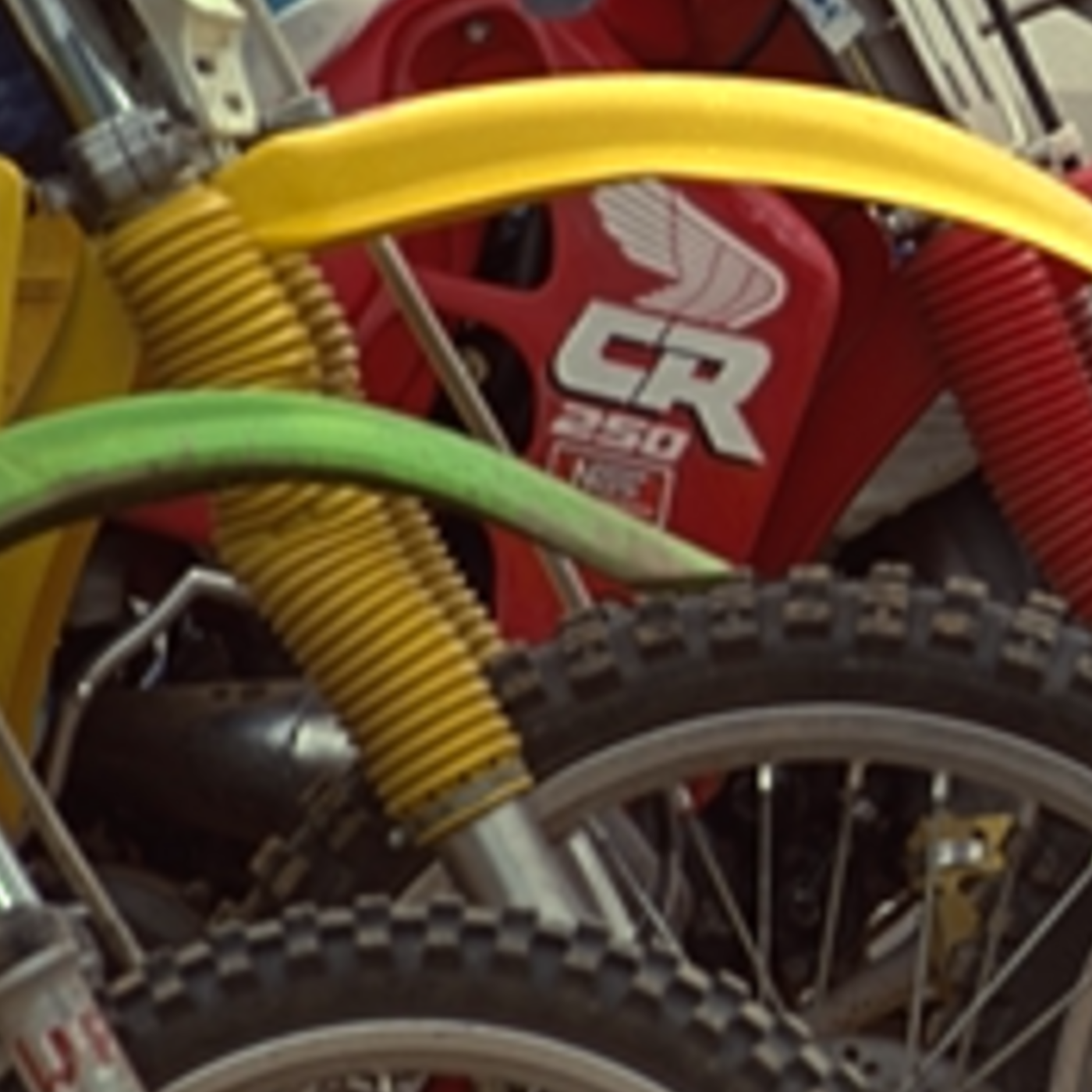}} \quad
			\subfigure[Mitchell]{\includegraphics[width=0.30\textwidth]{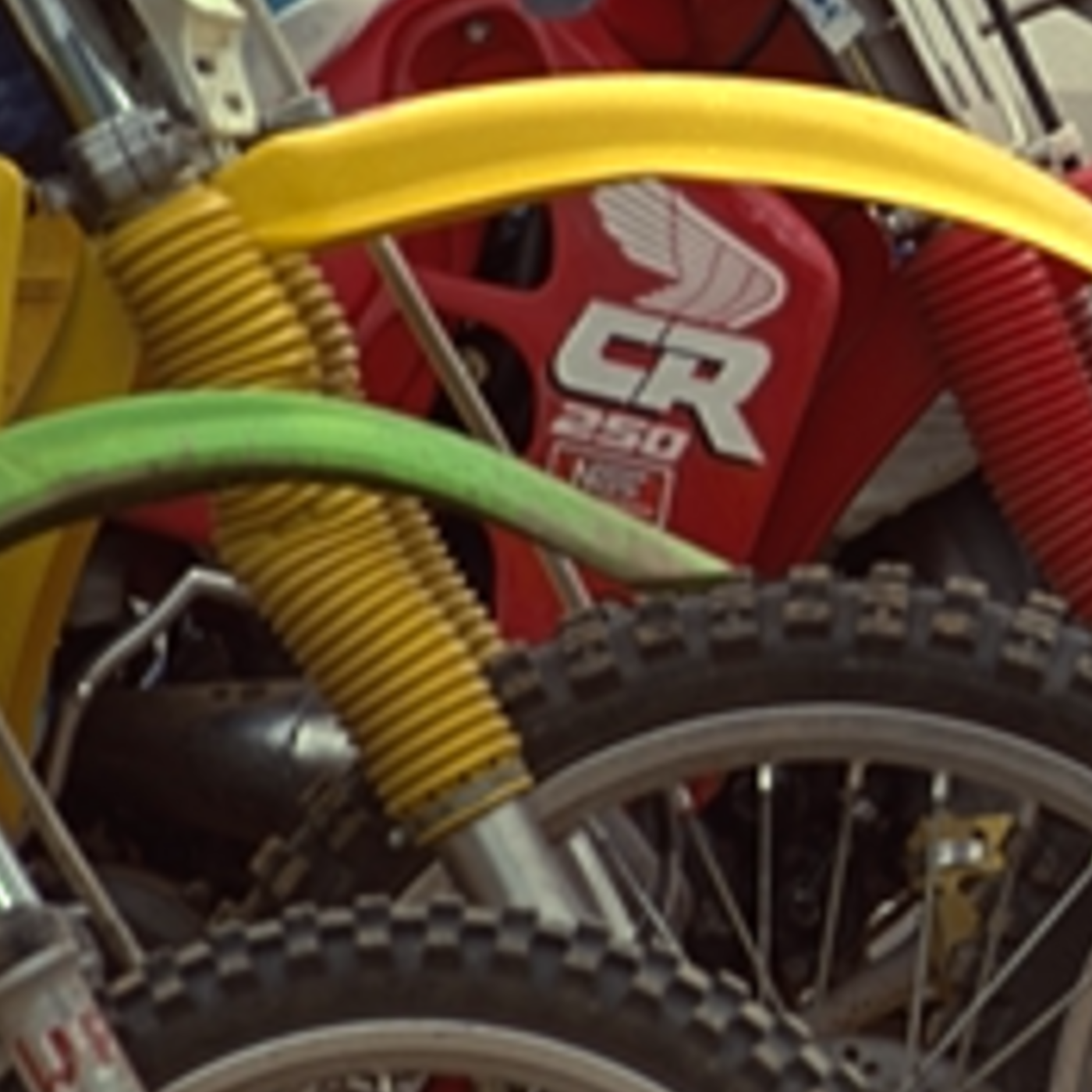}} \quad
			\subfigure[{WD} WENO]{\includegraphics[width=0.30\textwidth]{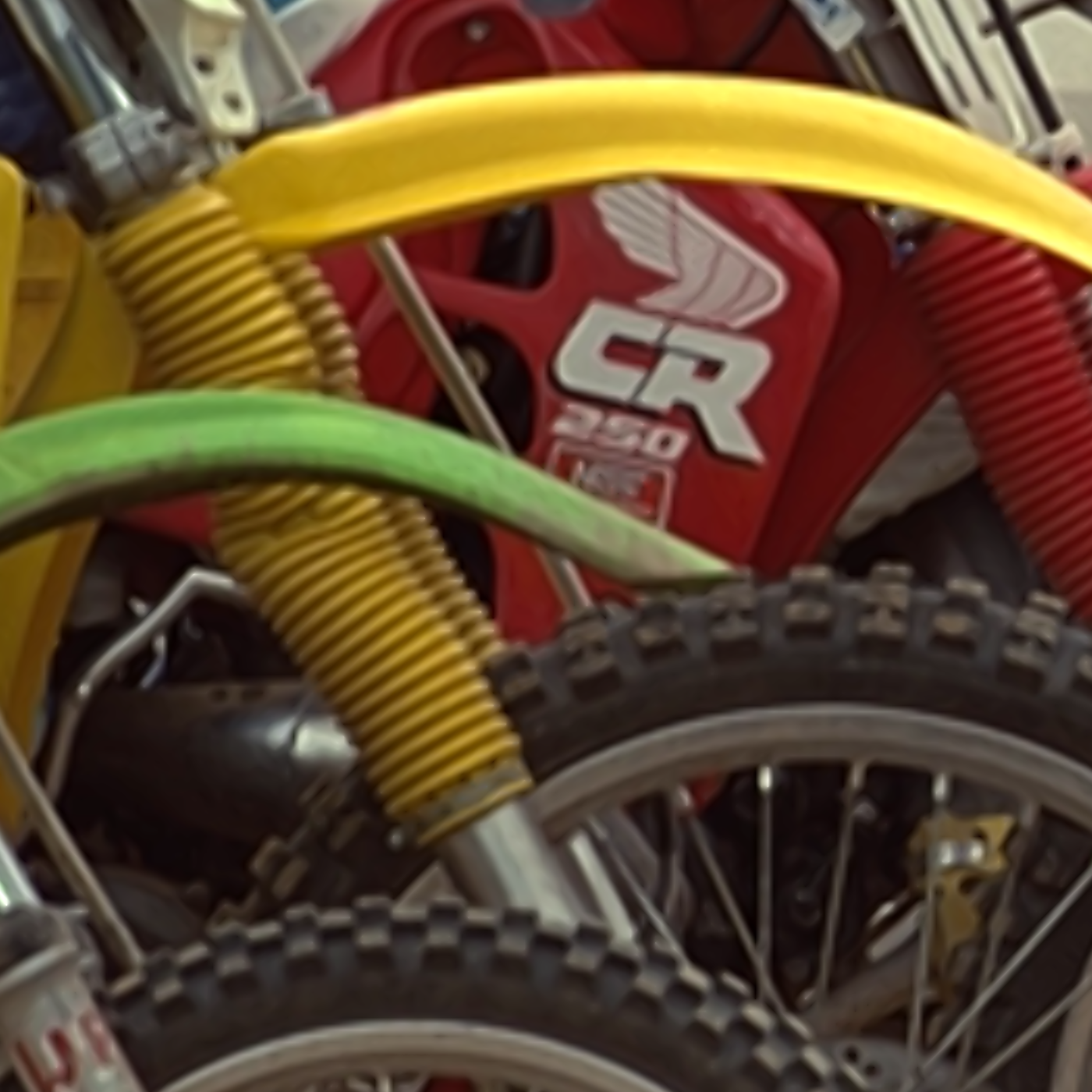}} \quad
			\subfigure[GCS]{\includegraphics[width=0.30\textwidth]{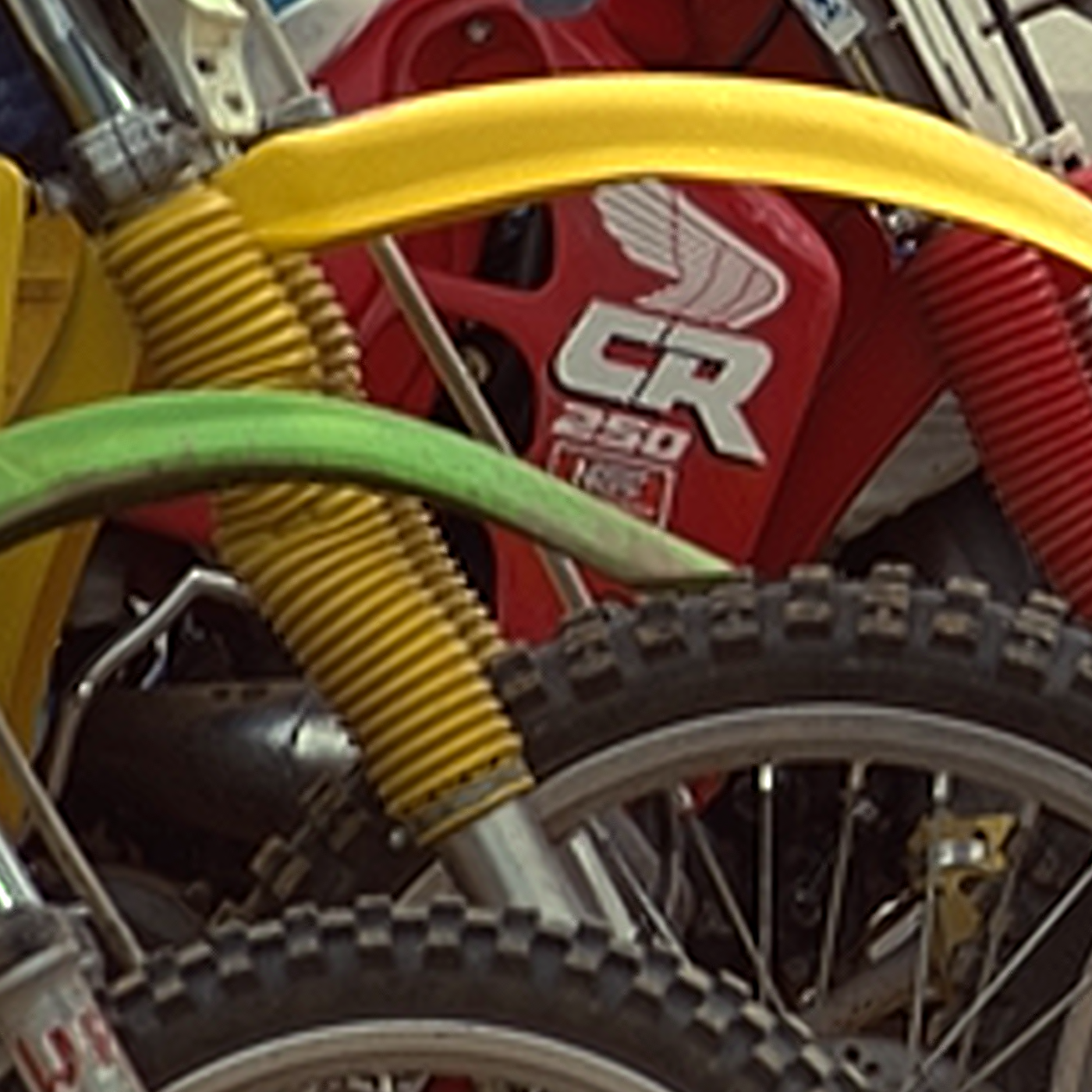}} \quad
			\caption{Kodak image No. 05, scale $d=8$}
			\label{fig:kodim05:scale8}
		\end{center}
	\end{figure}

\begin{figure}
		\begin{center}
			\subfigure[Input image]{\includegraphics[width=0.30\textwidth]{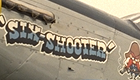}} \quad 
			\subfigure[Catrom]{\includegraphics[width=0.30\textwidth]{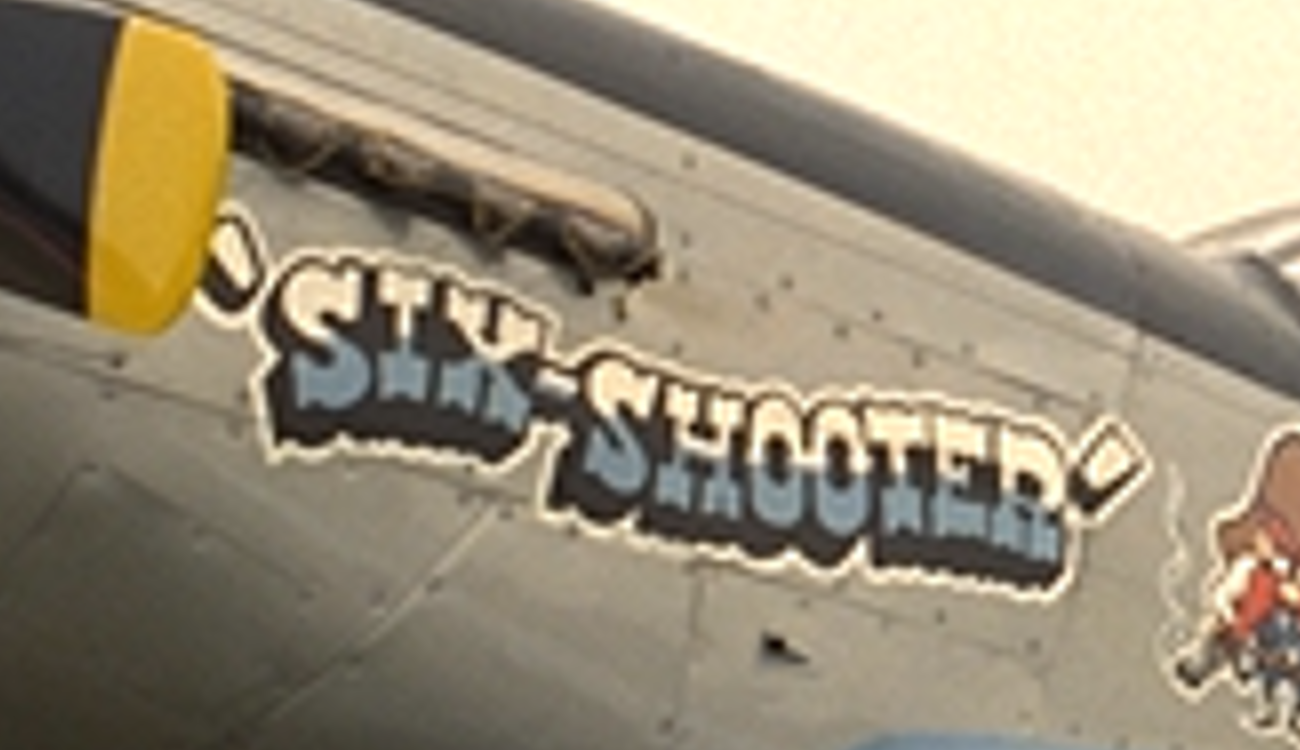}} \quad 
			\subfigure[Robidoux]{\includegraphics[width=0.30\textwidth]{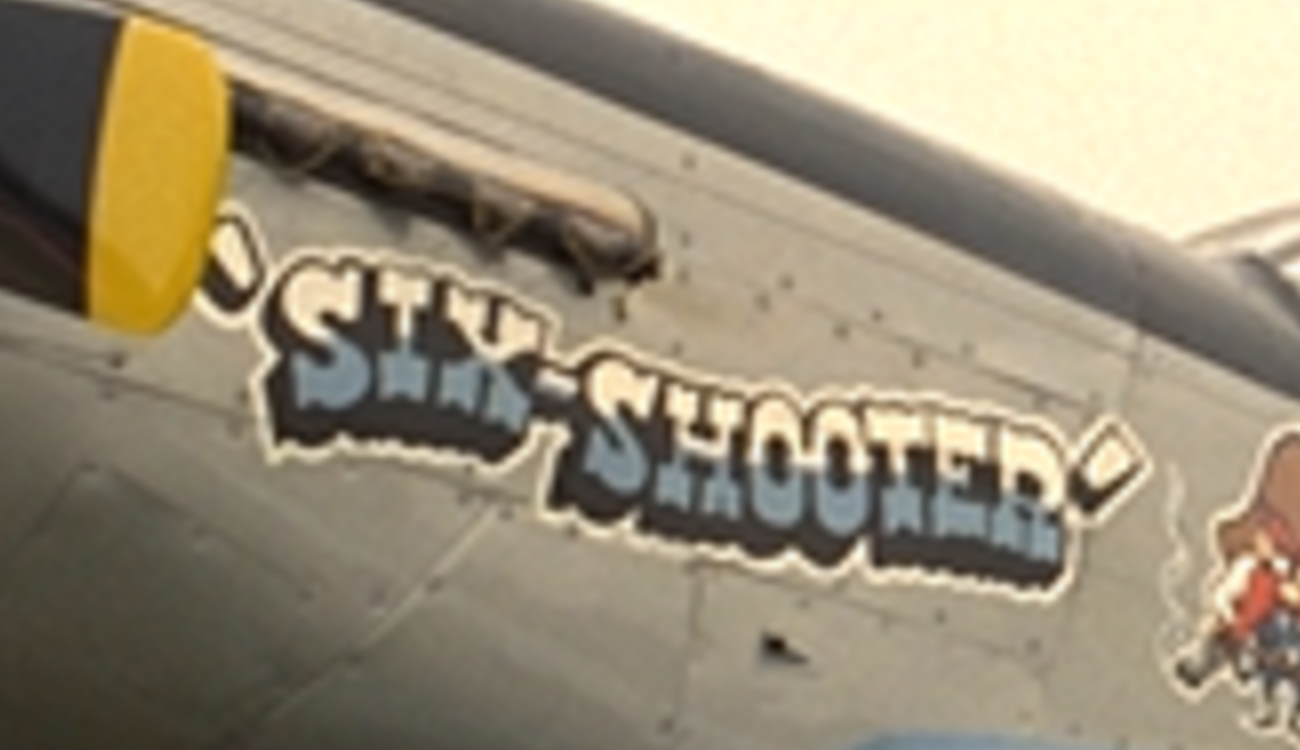}} \quad
			\subfigure[Mitchell]{\includegraphics[width=0.30\textwidth]{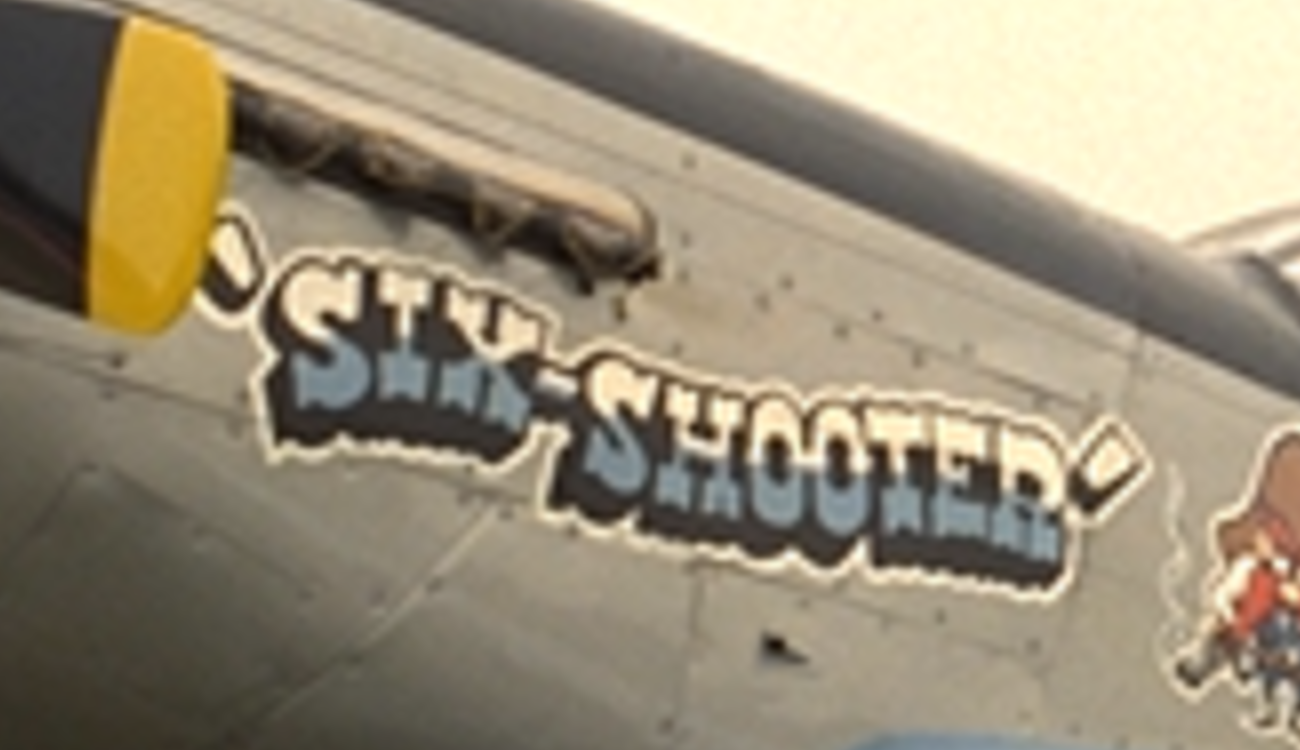}} \quad
			\subfigure[{WD} WENO]{\includegraphics[width=0.30\textwidth]{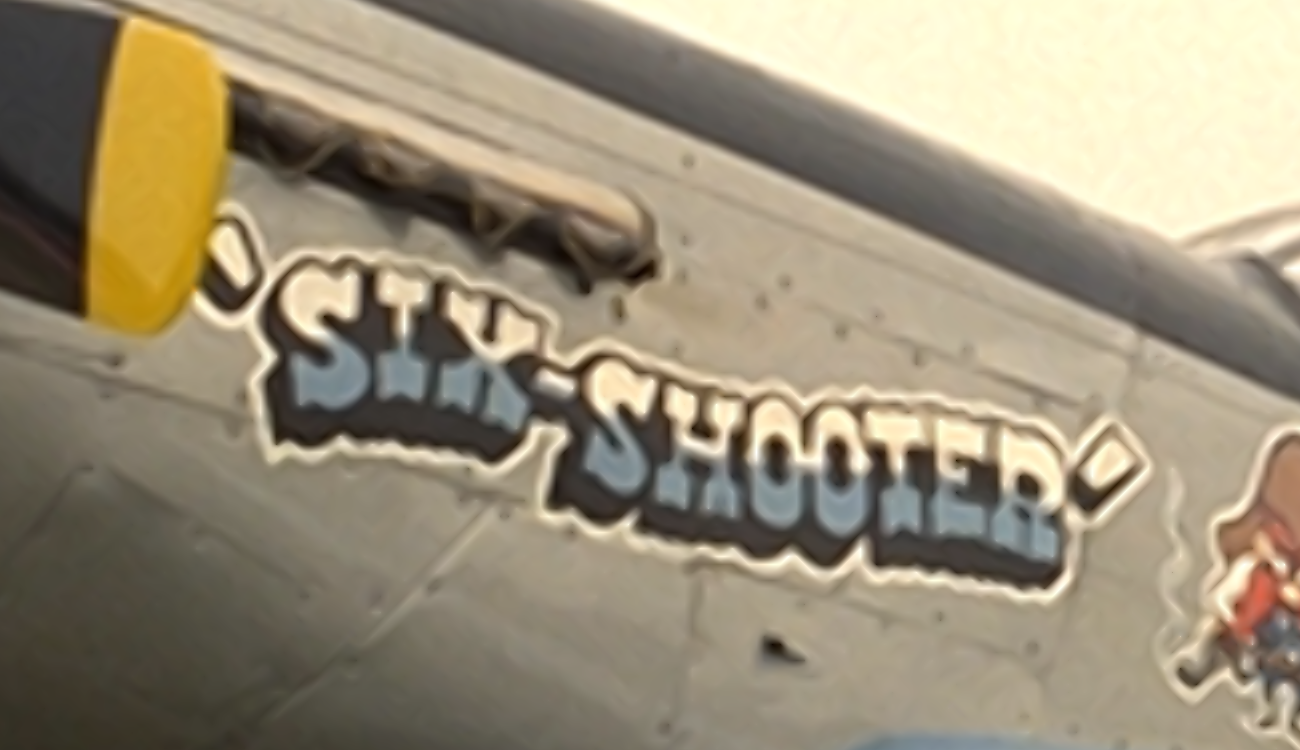}} \quad
			\subfigure[GCS]{\includegraphics[width=0.30\textwidth]{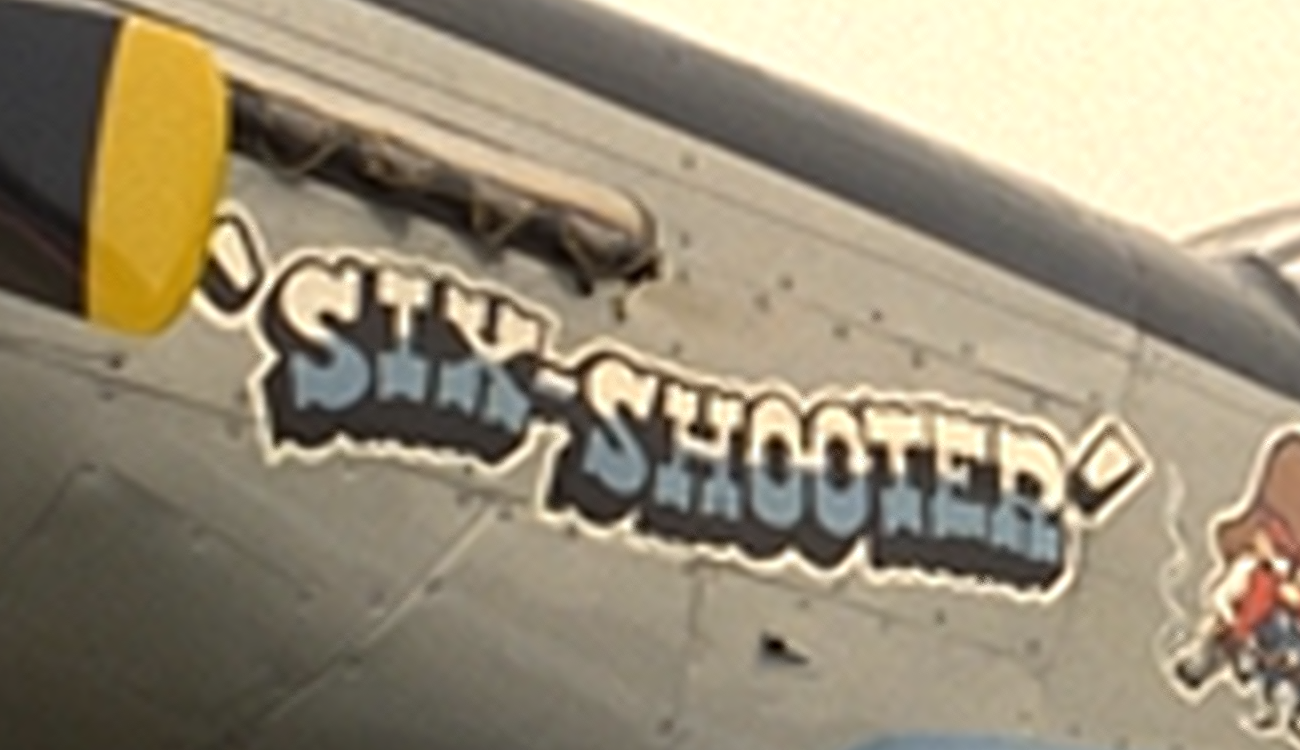}} \quad
			\caption{Kodak image No. 20, scale $d=8$}
			\label{fig:kodim20:scale8}
		\end{center}
	\end{figure}

Visual evaluation of the performance of {WD} WENO was performed by upsampling the original Kodak images No. $05$ and $20$ by a factor of $d = 8$, as shown in Figures \ref{fig:kodim05:scale8} and \ref{fig:kodim20:scale8}. Among all the methods presented, {WD} WENO and GCS stand out. The {WD} WENO images are very sharp, suggesting low diffusion, and have no detectable ringing effects. The GCS method, however, appears to produce the most visually sharp images, but the approximation has strong ringing effects that are noticeable in strong color transitions.
	
	\subsection{Images with sharp color transitions}
	
\begin{figure}
		\centering
		\subfigure[Pixel art]{\includegraphics[height=3cm]{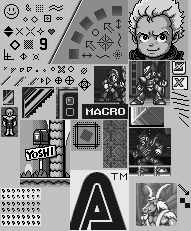}\label{fig:atari:input}}
		\subfigure[Tiger]{\includegraphics[height=3cm]{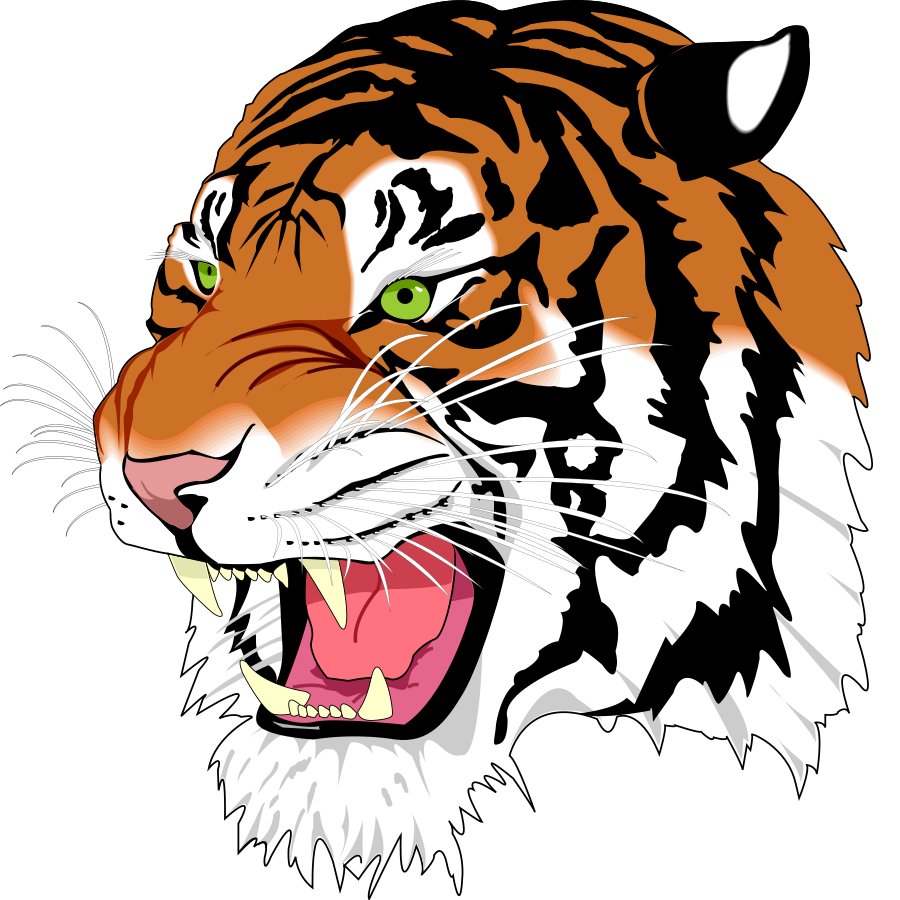}\label{fig:tiger:input}}
		\subfigure[Head]{\includegraphics[height=3cm]{./ZoomImages/fig_02_head/fig_02_head.png}\label{fig:head}}
		\caption{Selected test images}
		\label{fig:referenceimages}
	\end{figure}

	We selected two test images (Figures \ref{fig:atari:input} and \ref{fig:tiger:input}) to test the ability of the proposed {WD} WENO method to avoid oscillations and ringing effects and to produce reasonably sharp images.

	First, we chose the Tiger image in vector format\footnote{Available online at \texttt{https://commons.wikimedia.org/wiki/File:Ghostscript\_Tiger.svg}}. The vector image is then resampled to the following 
	resolution ($512$x$512$), which is chosen as the input image. Then, the vector image is resampled to 
	($1023$x$1023$), ($2045$x$2045$) and ($4089$x$4089$), which roughly corresponds with the scales of  
	$d=2, 4,$ and $8$, respectively. More precisely, resolutions are obtained with the following expression: 
	$$\mathrm{scaled}_{sz} = d(\mathrm{input}_{sz} - 1)+1.$$
			
	\begin{table}
		\begin{center}
			\caption{PSNR and MSSIM on the Tiger image for $d=4$ and $d=8$}
			\label{table:tiger:scale48}
			\begin{tabular}{lrrrr|rrrr}
				\hline
				& \multicolumn{4}{c}{$d=4$}& \multicolumn{4}{c}{$d=8$} \\
				Method          & PSNR      &Rank   &MSSIM      &Rank   &PSNR       &Rank   &MSSIM      &Rank  \\
				\hline
				{WD} WENO          &  21.6621  &  28  &  0.9001  &   9 &  21.3700  &  25  &  0.8932  &  2  \\
				\hline
				Tensor WENO   &  21.2307  &  33  &  0.8871  &  26 &  20.9963  &  29  &  0.8822  &  20  \\
				\hline
				Box           &  21.5205  &  29  &  0.8814  &  27 &  21.5205  &  30  &  0.8814  &  27 \\
				Triangle      &  22.5313  &  25  &  0.8875  &  25 &  21.5008  &  24  &  0.8797  &  25 \\
				\hline                                           
				Catrom        &  23.2599  &  17  &  0.9031  &   5 &  22.0893  &  17  &  0.8906  &  4  \\
				Cubic         &  21.4835  &  31  &  0.8611  &  31 &  20.6357  &  32  &  0.8618  &  31 \\
				Gaussian      &  22.0278  &  27  &  0.8755  &  29 &  21.0895  &  27  &  0.8712  &  29 \\
				Hermite       &  22.6928  &  22  &  0.8936  &  21 &  21.6377  &  22  &  0.8846  &  13 \\
				Lagrange      &  23.1486  &  19  &  0.9004  &   9 &  21.9990  &  18  &  0.8889  &  7  \\
				Mitchell      &  22.7068  &  21  &  0.8913  &  22 &  21.6446  &  21  &  0.8825  &  19 \\
				Quadratic     &  21.9728  &  28  &  0.8738  &  30 &  21.0441  &  28  &  0.8704  &  30 \\
				Robidoux      &  22.6275  &  24  &  0.8895  &  24 &  21.5803  &  23  &  0.8812  &  23 \\
				RobidouxSharp &  22.8310  &  20  &  0.8941  &  19 &  21.7449  &  19  &  0.8844  &  14 \\
				Spline        &  21.4835  &  31  &  0.8611  &  31 &  20.6357  &  32  &  0.8618  &  31 \\
				\hline
				Bartlett      &  23.5207  &  11  &  0.8945  &  15 &  22.2953  &  11  &  0.8818  &  21 \\
				Blackman      &  23.4931  &  12  &  0.9013  &   8 &  22.2734  &  12  &  0.8875  &  9  \\
				Bohman        &  23.4826  &  13  &  0.9018  &   7 &  22.2651  &  13  &  0.8880  &  8  \\
				\hline
				Cosine        &  23.5930  &   6  &  0.8941  &  19 &  22.3525  &   6  &  0.8807  &  24 \\
				Hamming       &  23.5992  &   4  &  0.8967  &  14 &  22.3573  &   4  &  0.8834  &  16 \\
				Hann          &  23.5535  &   9  &  0.8944  &  16 &  22.3211  &   9  &  0.8813  &  22 \\
				Jinc          &  22.3566  &  26  &  0.8607  &  33 &  21.3597  &  26  &  0.8554  &  33 \\
				Kaiser        &  23.5449  &  10  &  0.8978  &  13 &  22.3142  &  10  &  0.8843  &  15 \\
				Lanczos       &  23.5717  &   7  &  0.8989  &  11 &  22.3355  &   7  &  0.8850  &  11 \\
				Lanczos2      &  23.2776  &  16  &  0.9025  &   6 &  22.1031  &  16  &  0.8897  &  6  \\
				Lanczos2Sharp &  23.3466  &  15  &  0.9044  &   2 &  22.1594  &  15  &  0.8910  &  3  \\
				LanczosRadius &  23.5717  &   7  &  0.8989  &  11 &  22.3355  &   7  &  0.8850  &  11 \\
				LanczosSharp  &  23.6128  &   3  &  0.8997  &  10 &  22.3680  &   3  &  0.8852  &  10 \\
				Parzen        &  23.4421  &  14  &  0.9035  &   4 &  22.2333  &  14  &  0.8899  &  5  \\
				Point         &  21.5205  &  29  &  0.8814  &  27 &  20.7240  &  30  &  0.8758  &  27 \\
				Sinc          &  23.9279  &   1  &  0.8944  &  16 &  22.6267  &   1  &  0.8834  &  16 \\
				SincFast      &  23.9279  &   1  &  0.8944  &  16 &  22.6267  &   1  &  0.8834  &  16 \\
				Welch         &  23.5978  &   5  &  0.8905  &  23 &  22.3565  &   5  &  0.8777  &  26 \\
				\hline
				GCS           &  23.2294  &  18  &  0.9083  &   1 &  21.7315  &  20  &  0.8938  &  1  \\
				\hline
			\end{tabular}
		\end{center}
	\end{table}	
	
Table \ref{table:tiger:scale48} shows the PSNR and MSSIM results for the Tiger image for two scaling factors.
It is worth noting that the results for $d=2$ have been omitted because they do not show significant differences in PSNR or MSSIM ranks.
The images obtained with the {WD} WENO method do not give very good results in the PSNR norm, while the opposite is true for the same results in the MSSIM norm.
However, the obtained images are quite sharp and do not show ringing effects.

\begin{figure}
		\begin{center}
			\subfigure[Input image]{\includegraphics[width=0.22\textwidth]{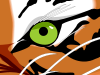}} \quad 
			\subfigure[Catrom]{\includegraphics[width=0.22\textwidth]{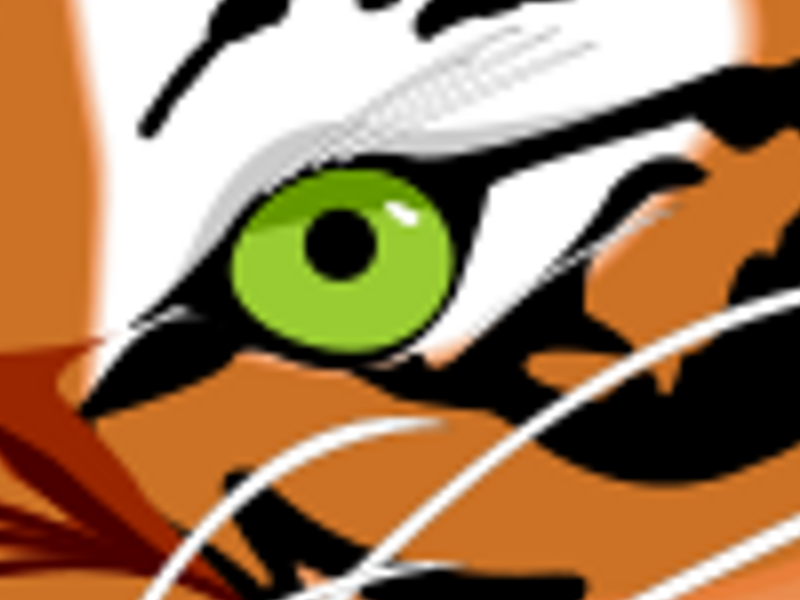}} \quad 
			\subfigure[Parzen]{\includegraphics[width=0.22\textwidth]{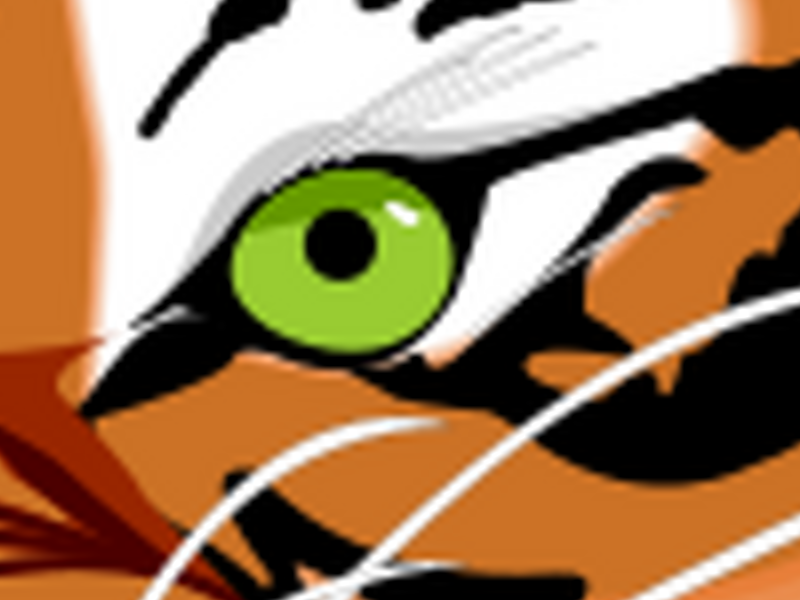}} \quad
			\subfigure[Lanczos2Sharp]{\includegraphics[width=0.22\textwidth]{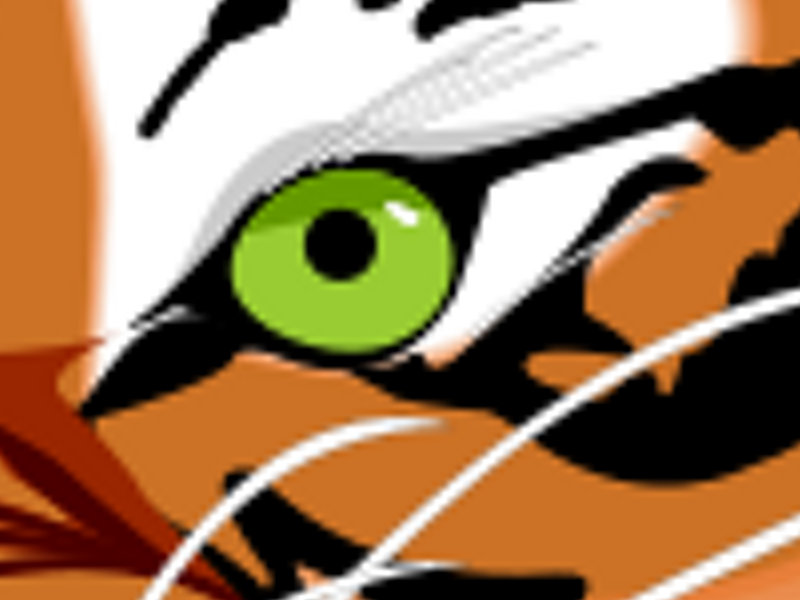}} \quad
			\subfigure[Spline]{\includegraphics[width=0.22\textwidth]{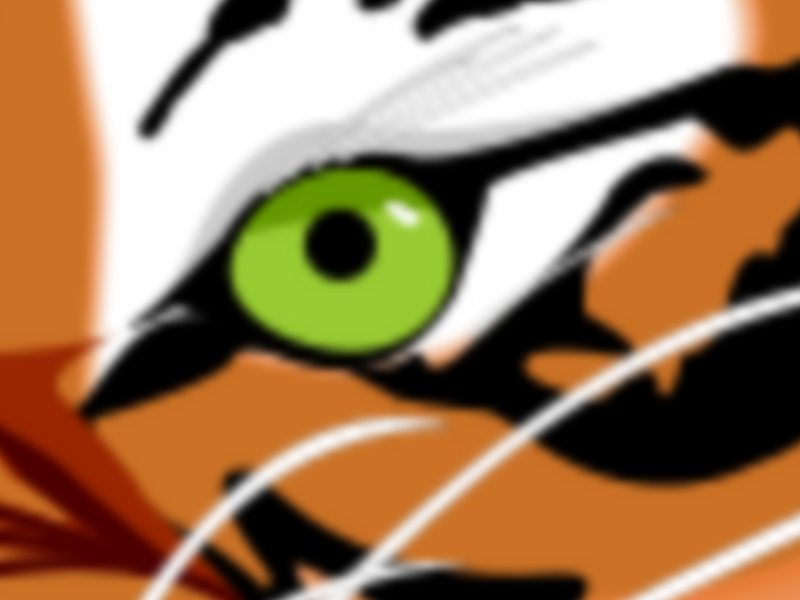}} \quad
			\subfigure[Tensor WENO]{\includegraphics[width=0.22\textwidth]{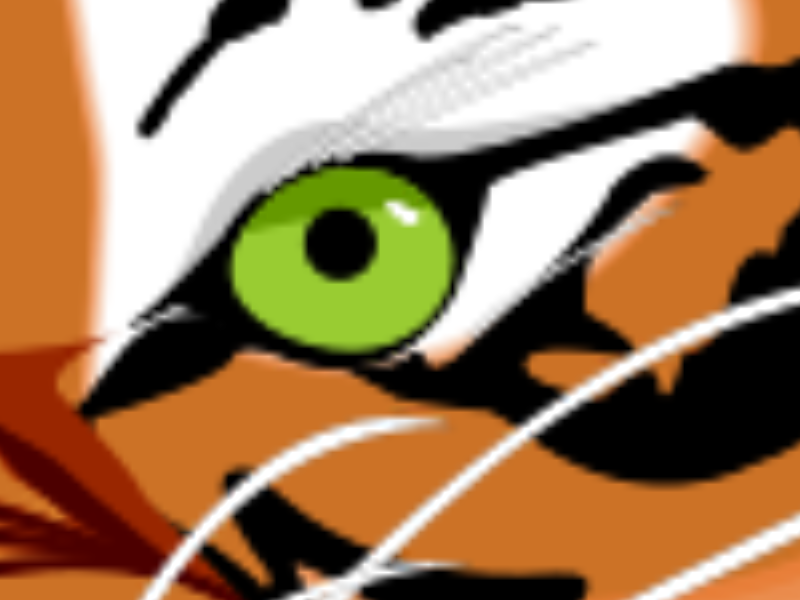}} \quad
			\subfigure[{WD} WENO]{\includegraphics[width=0.22\textwidth]{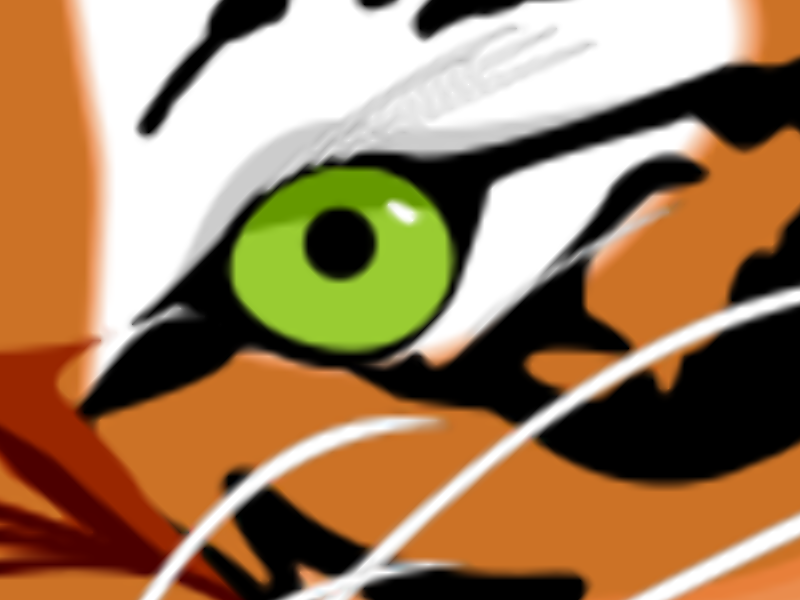}} \quad
			\subfigure[GCS]{\includegraphics[width=0.22\textwidth]{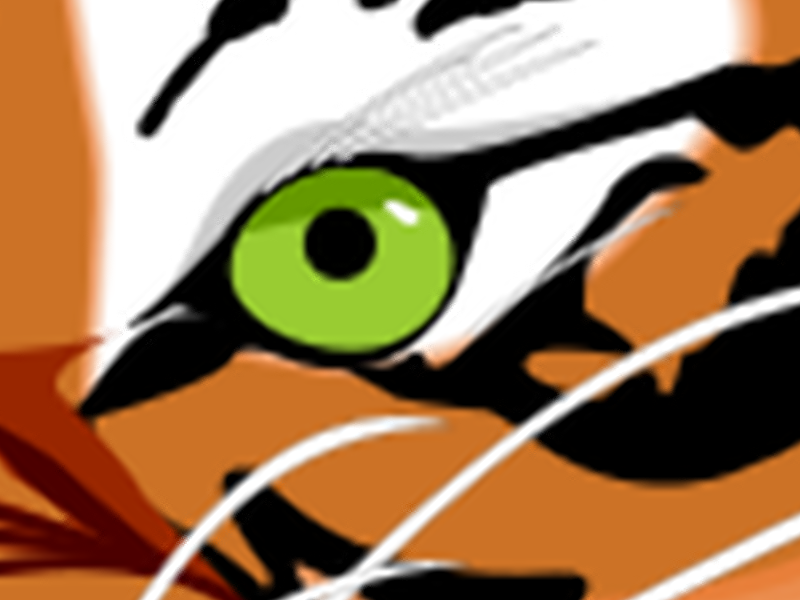}} \quad
			
			\caption{Eye of the tiger, scale $d=8$}
			\label{fig:tiger:scale8}
		\end{center}
	\end{figure}

Figure \ref{fig:tiger:scale8} shows the comparison of the different methods for scale $d=8$.
Although the GCS method gives the best results compared to the other methods, {WD} WENO introduces the least amount of diffusion and does not produce ringing effects.
The ringing effect is also visible with the GCS method. The same effect is also observed with other scaling factors and is therefore not shown. {We should point out that WD WENO produces significantly better results quantitatively and visually than the tensor WENO and the example of this comparison can be seen in Figure \ref{fig:tiger:scale8} and Table \ref{table:tiger:scale48}.} 
	
\begin{figure}
		\begin{center}
			\subfigure[Input image]{\includegraphics[width=0.30\textwidth]{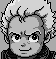}} \quad 
			\subfigure[Catrom]{\includegraphics[width=0.30\textwidth]{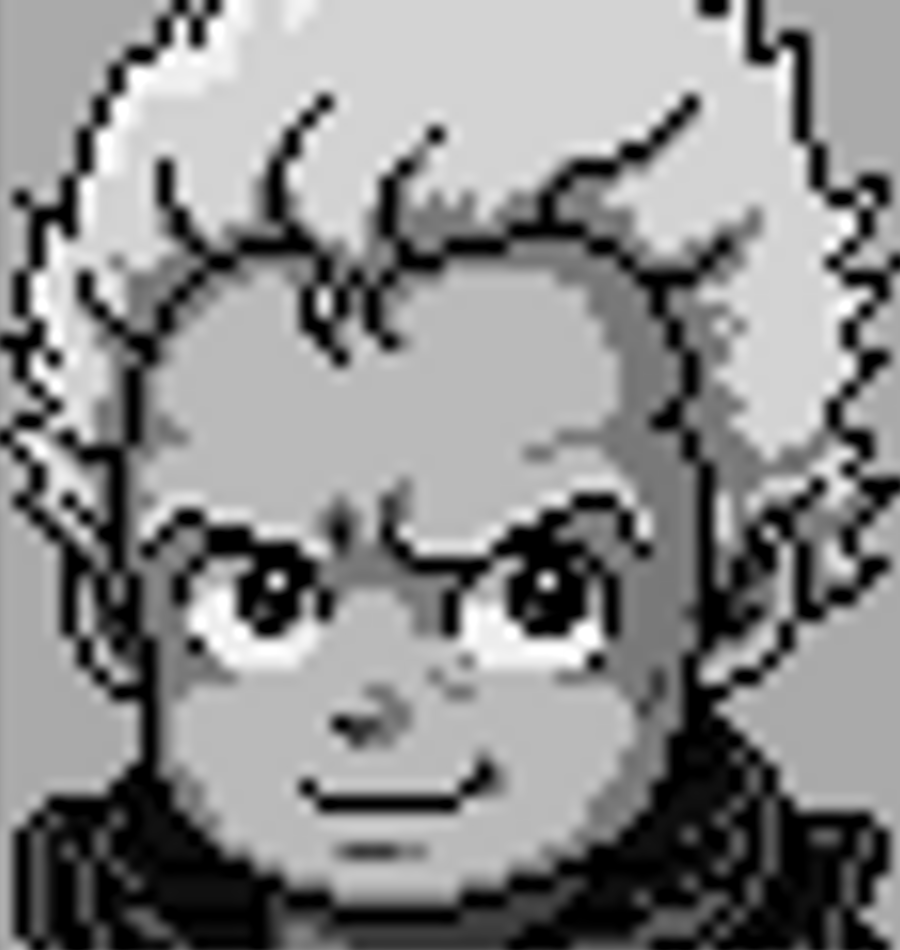}} \quad 
			\subfigure[Lanczos Radius]{\includegraphics[width=0.30\textwidth]{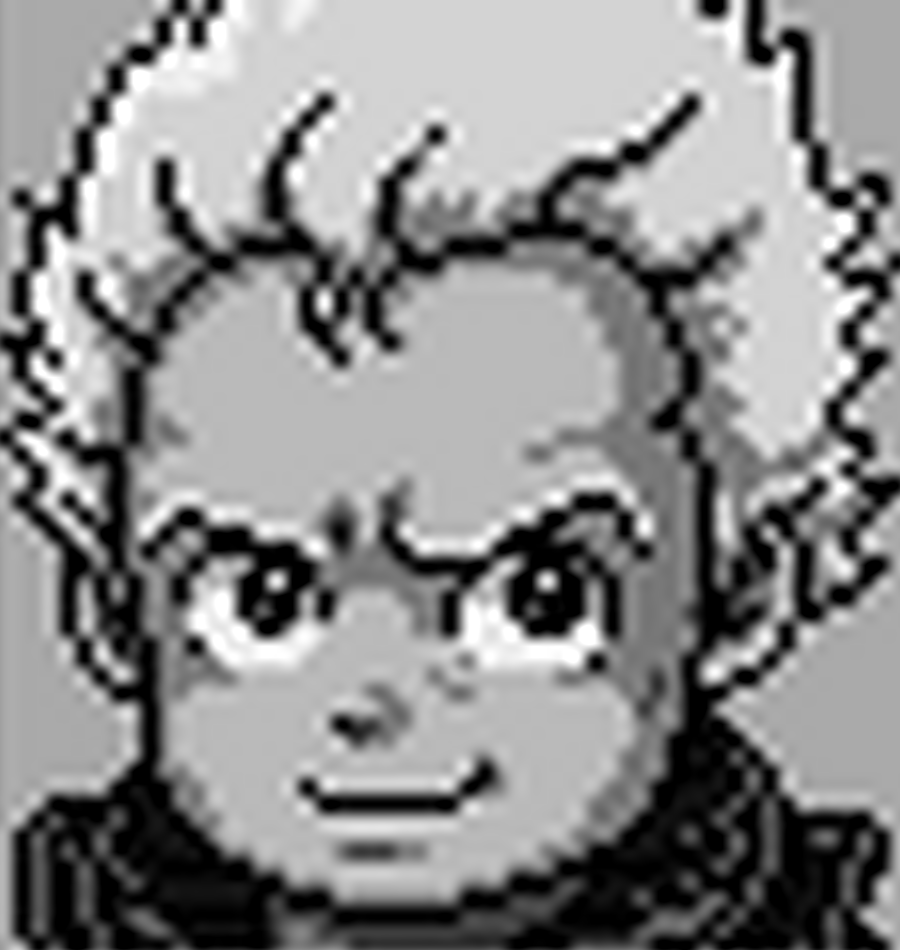}} \quad
			\subfigure[Spline]{\includegraphics[width=0.30\textwidth]{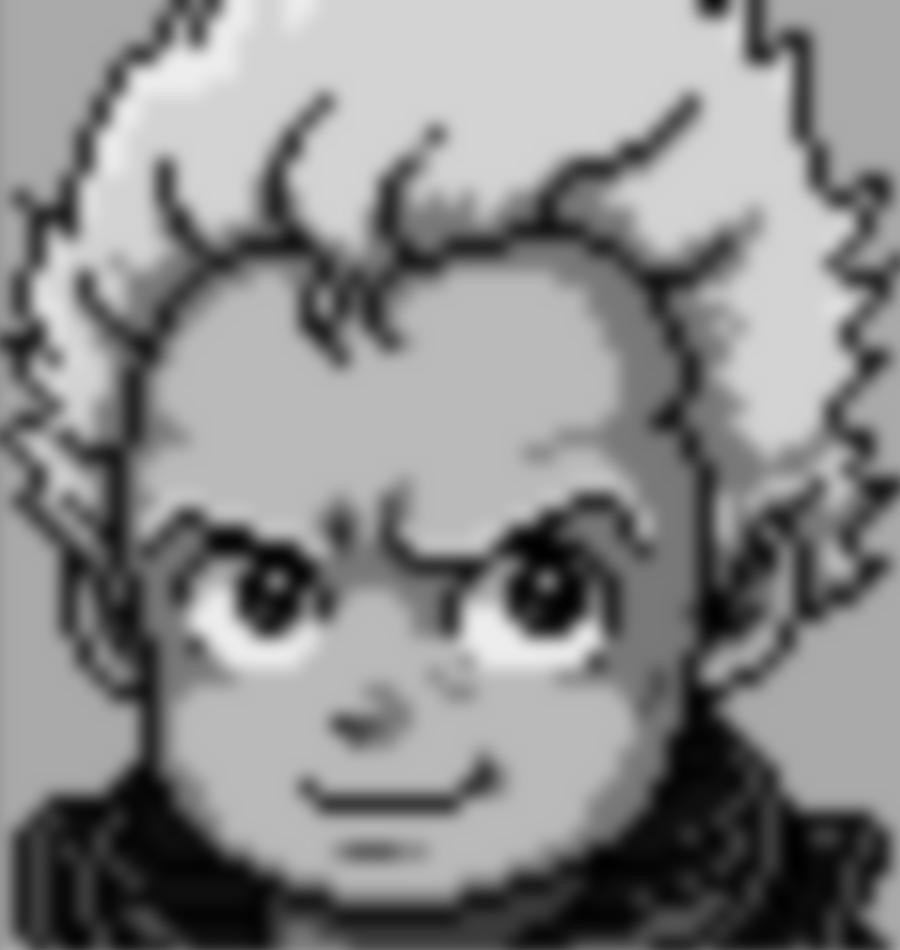}} \quad
			\subfigure[{WD} WENO]{\includegraphics[width=0.30\textwidth]{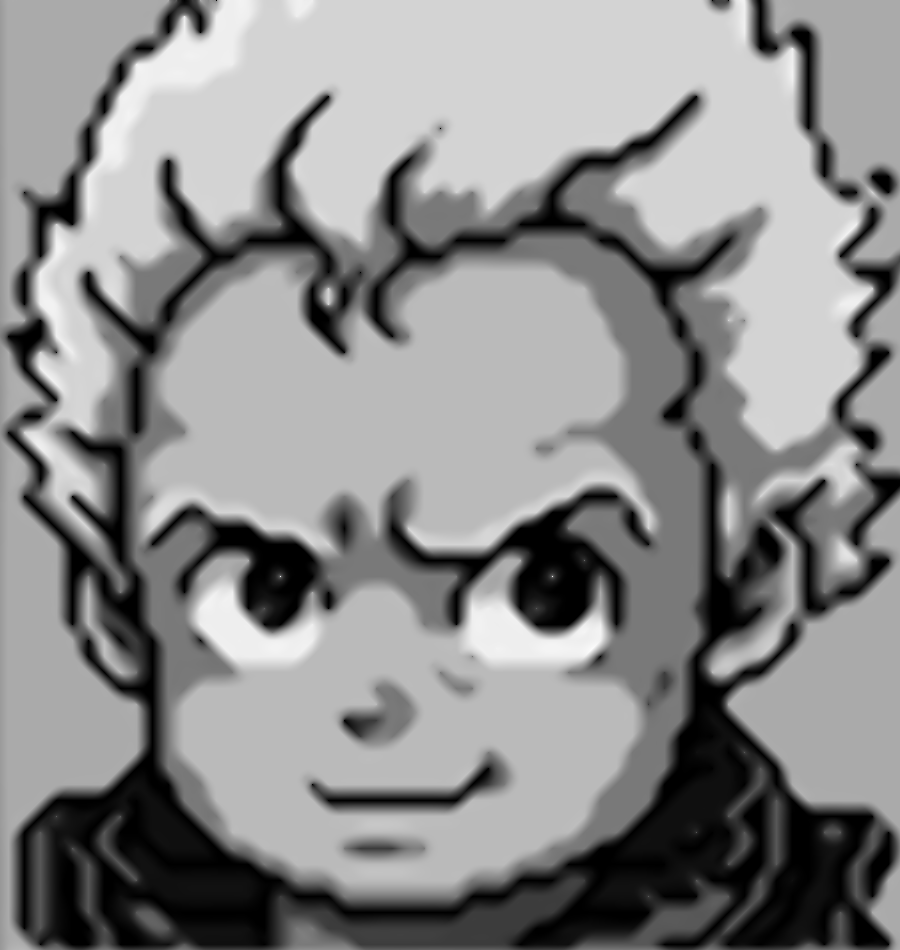}} \quad
			\subfigure[GCS]{\includegraphics[width=0.30\textwidth]{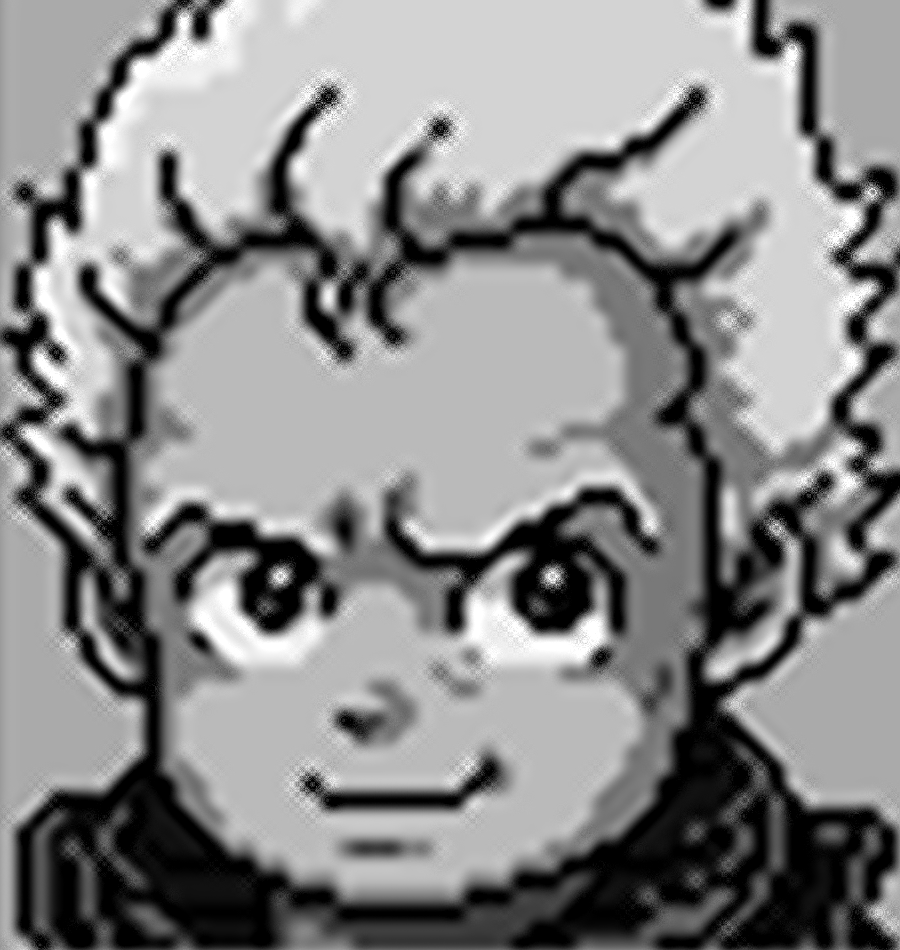}} \quad
			\caption{Pixel art, Head, scale $d=16$}
			\label{fig:yoshi:scale16}
		\end{center}
	\end{figure}

\begin{figure}
		\begin{center}
			\subfigure[Input image]{\includegraphics[width=0.30\textwidth]{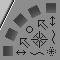}} \quad 
			\subfigure[Catrom]{\includegraphics[width=0.30\textwidth]{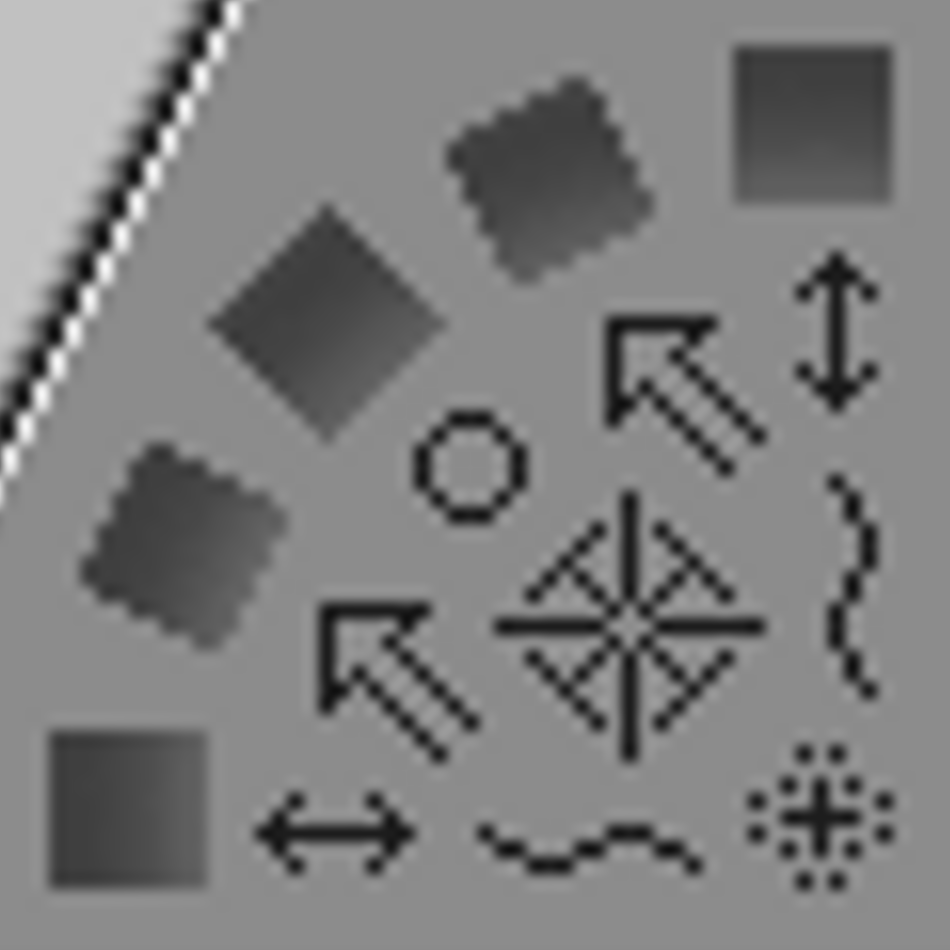}} \quad 
			\subfigure[Lanczos Radius]{\includegraphics[width=0.30\textwidth]{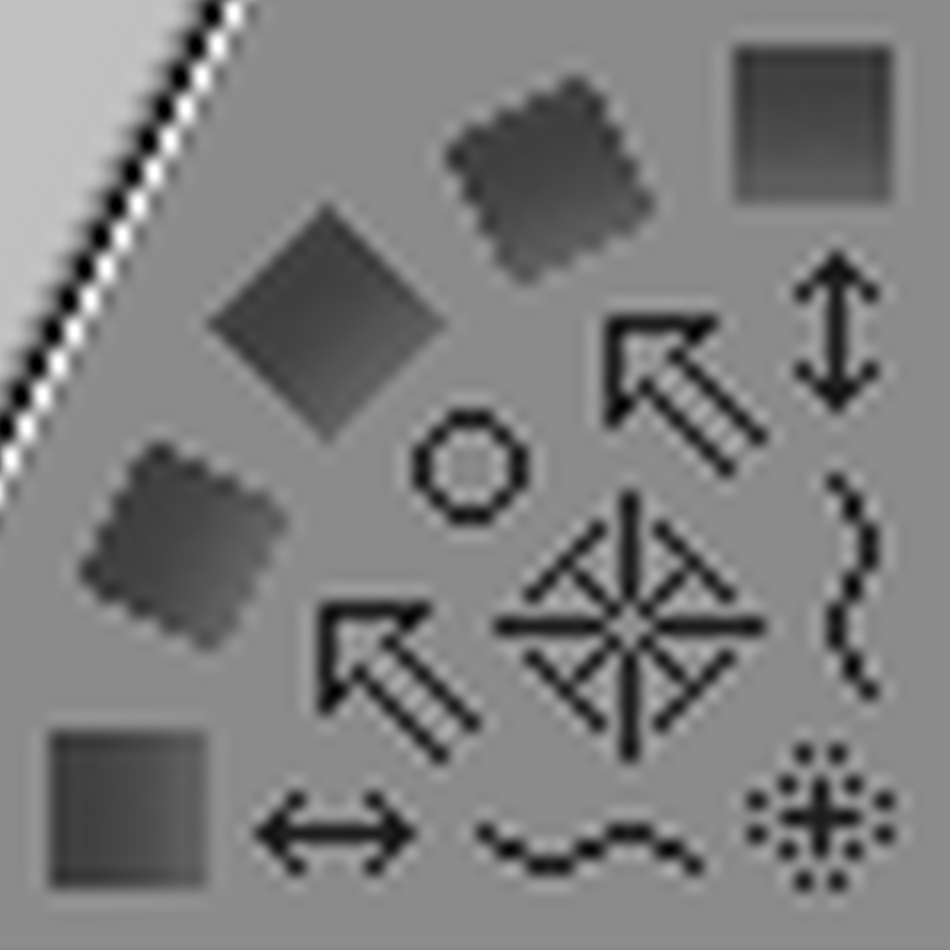}} \quad
			\subfigure[Spline]{\includegraphics[width=0.30\textwidth]{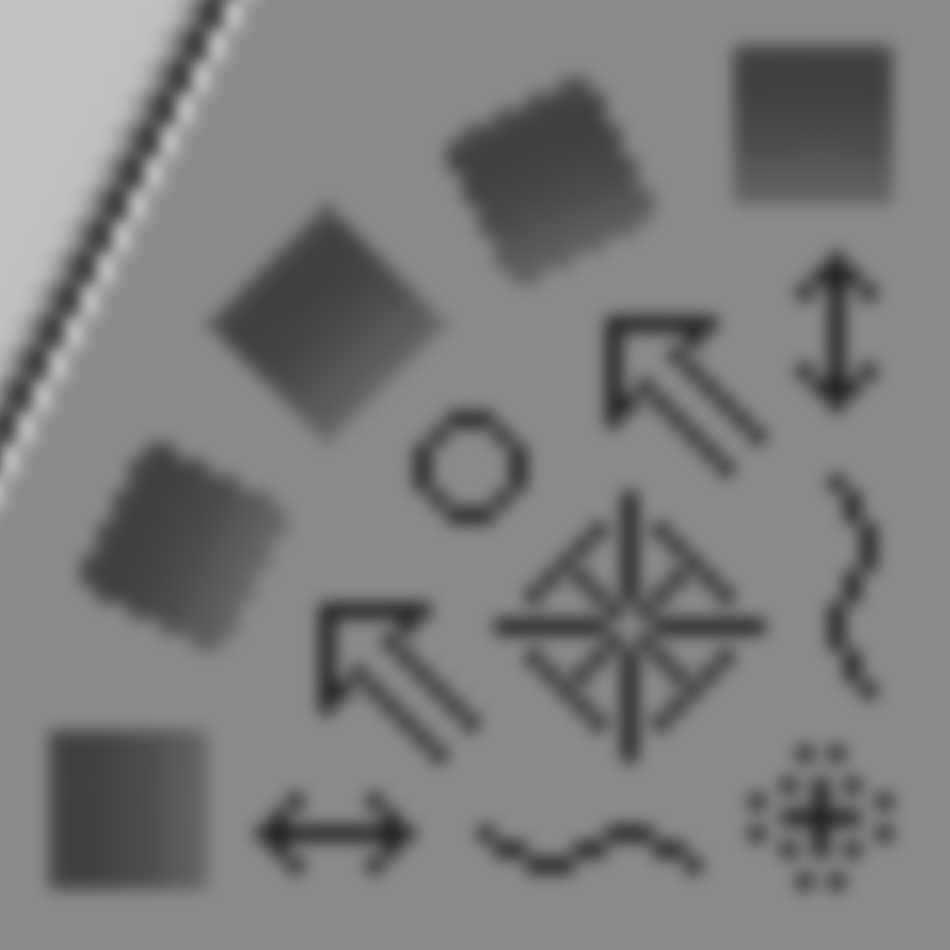}} \quad
			\subfigure[{WD} WENO]{\includegraphics[width=0.30\textwidth]{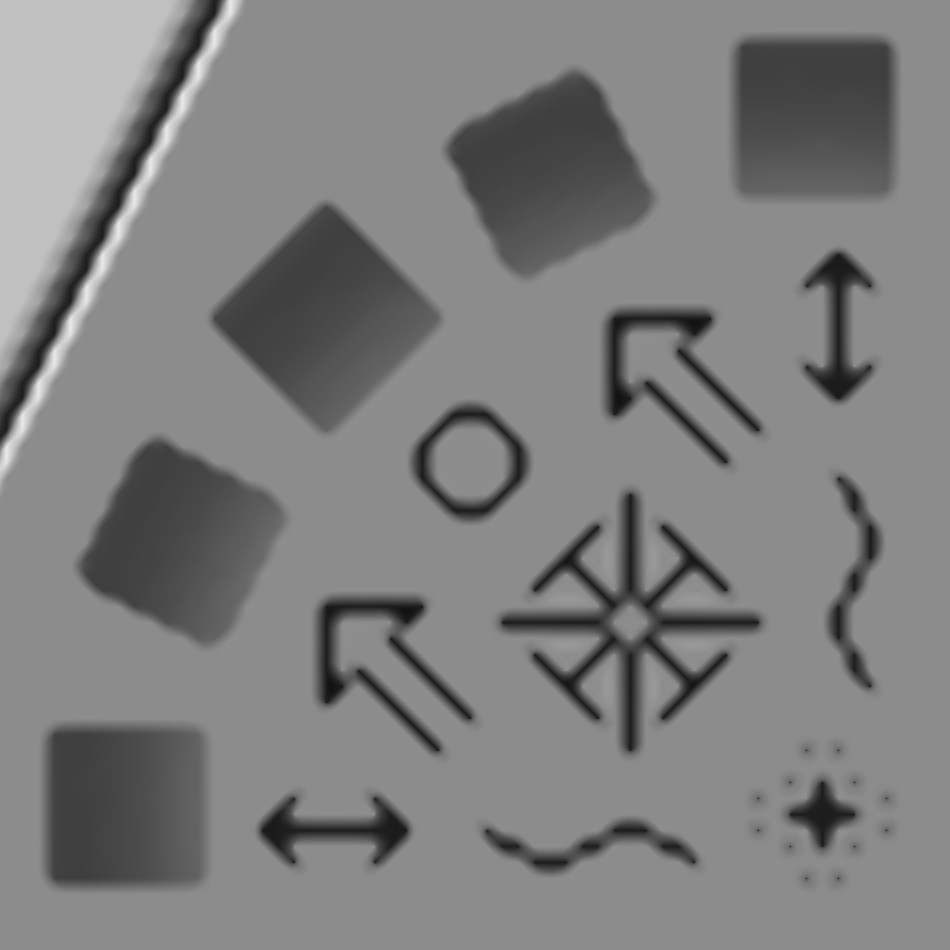}} \quad
			\subfigure[GCS]{\includegraphics[width=0.30\textwidth]{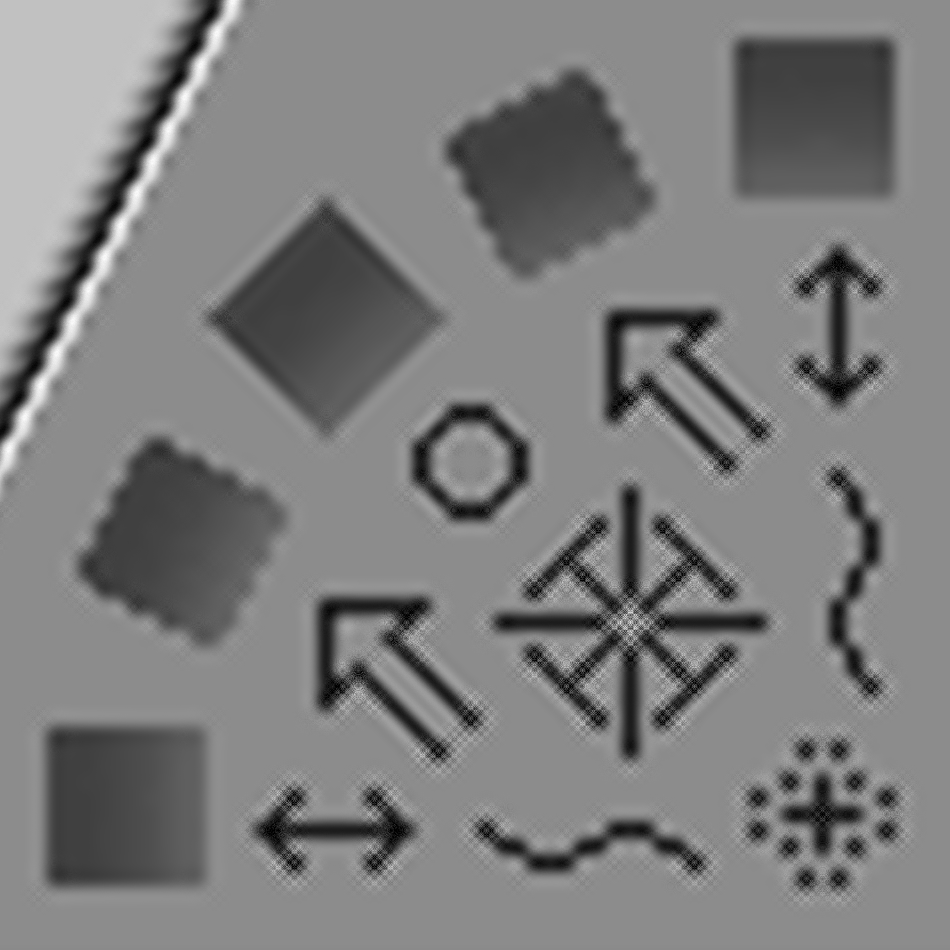}} \quad
			\caption{Pixel art, squares, scale $d=16$}
			\label{fig:squares:scale16}
		\end{center}
	\end{figure}

To further evaluate the performance of the {WD} WENO method, the pixel art image shown in Figure \ref{fig:atari:input} was scaled for $d=16$.
Figures \ref{fig:yoshi:scale16} and \ref{fig:squares:scale16} clearly show that the {WD} WENO method has the best overall performance. {WD} WENO reduces the staircase effect significantly for multiple angles of lines which is far better then expected. It is interesting to note that GCS method quite fails with the image from the Figure \ref{fig:atari:input}, since next to the ringing effect, it shows staircase effect on the slanted squares, which was not expected.  
	
	\section{Conclusion}
	\label{sec:conclusion}
The aim of this work was to develop an efficient algorithm with low memory load that produces sharp images, with the method having fourth-order accuracy in smooth regions while maintaining a high degree of accuracy near discontinuities and ensuring good approximation properties when applied to digital images. The choice of a relatively low degree of accuracy, for WENO methods, and a small number of interpolation directions should make this method {with} comparable complexity to the tensor WENO method {and commonly used linear methods} while improving the quality of the interpolated digital images. For image doubling this goal was achieved, but for general resolutions it is at most twice as complex as {tensor WENO method}.

Tests with standard test images have confirmed that the algorithm does not produce excessive numerical artifacts such as ringing or aliasing effects and that the image quality and visual sharpness are almost equivalent to the much more complex GCS method.
For images with sharp color transitions, the {WD} WENO method achieved very good results compared to other methods in terms of standard metrics and outperformed all in visual qualitative inspection. The algorithm allows direct parallelization and vectorization of the algorithm and can be generalized for higher dimensions, multiple directions and higher accuracy, which will be the topic of future work.

\section*{Acknowledgements}
This research is primarily supported by the Croatian Science Foundation under the project IP-2019-04-1239.
\section*{Data availability}
The Python code needed to reproduce this research is available upon request.

	\bibliographystyle{spmpsci} 
	\bibliography{references}
\end{document}